# Subfactors of index less than 5, part 2: triple points

SCOTT MORRISON
DAVID PENNEYS
EMILY PETERS
NOAH SNYDER

URLs: http://tqft.net/ http://math.berkeley.edu/~dpenneys
http://math.mit.edu/~eep and http://math.columbia.edu/~nsnyder

Email: scott@tqft.net dpenneys@math.berkeley.edu eep@math.mit.edu nsnyder@math.columbia.edu

**Abstract** We summarize the known obstructions to subfactors with principal graphs which begin with a triple point. One is based on Jones's quadratic tangles techniques, although we apply it in a novel way. The other two are based on connections techniques; one due to Ocneanu, and the other previously unpublished, although likely known to Haagerup.

We then apply these obstructions to the classification of subfactors with index below 5. In particular, we eliminate three of the five families of possible principal graphs called "weeds" in the classification from [MS10].

## Contents

# 1 Introduction

This is the second paper in a project where we extend the classification of subfactors of small index. The first result about subfactors of small index is Jones's index theorem for subfactors [Jon83] which states that the index of a subfactor lies in the range $\{4\cos^2(\frac{\pi}{n})|n = 3, 4, \ldots\} \cup [4, \infty]$. Any of these numbers can be realized as the index of a subfactor whose standard invariant is either a quotient of Temperley-Lieb (if the index is less than $4$) [Jon83] or Temperley-Lieb (if the index is at least $4$) [Pop93]. However, once you ignore the subfactors with Temperley-Lieb standard invariant, the possible indices for irreducible subfactors are again quantized in an interval above $4$. Haagerup began the classification of subfactors with index 'only a little larger' than four in [Haa94]. In that paper, he showed there are no extremal subfactors (other than those with Temperley-Lieb standard invariant) with index in the range $(4, \frac{5+\sqrt{13}}{2})$. Furthermore, he gave a complete list of possible principal graphs of extremal subfactors whose index falls in the range $(4, 3+\sqrt{3})$. (He states the result up to $3+\sqrt{3}$, and proves it up to $3+\sqrt{2}$.) Most of the graphs on this list were excluded by Bisch [Bis98] and Asaeda-Yasuda [AY09], while the remaining $3$ graphs were shown to come from (unique hyperfinite) subfactors by Asaeda-Haagerup [AH99] and Bigelow-Morrison-Peters-Snyder [BMPS09]. Haagerup's classification stops at index $3+\sqrt{3}$ for reasons of computational convenience, and because a Goodman-de la Harpe-Jones subfactor [GdlHJ89] was already known to exist at that index. In this series of papers we extend this classification to index $5$.

**Theorem 1.1** *There are exactly ten subfactor planar algebras other than Temperley-Lieb with index between $4$ and $5$: the Haagerup planar algebra and its dual [AH99], the extended Haagerup planar algebra and its dual [BMPS09], the Asaeda-Haagerup planar algebra [AH99] and its dual, the 3311 Goodman-de la Harpe-Jones planar algebra [GdlHJ89], and Izumi's self-dual 2221 planar algebra [Izu01] and its complex conjugate.*

In the first paper of this series [MS10], Morrison and Snyder gave an initial classification result analogous to Haagerup's classification. In order to complete the classification we rule out the remaining families in the subsequent papers. Thus parts 2 through 4 are closer in sprit to the papers of Bisch [Bis98] and Asaeda–Yasuda [AY09] which eliminated certain families of candidate principal graphs coming from Haagerup's classification.

In order to state our main results we rapidly recall some terminology from part 1. A translation of a graph pair is used to indicate a graph pair obtained by increasing the supertransitivity by an even integer (the supertransitivity is the number of edges between the initial vertex and the first vertex of degree more than two). An extension of a graph pair is a graph pair obtained by extending the graphs in any way at greater depths (i.e. adding vertices and edges at the right), even infinitely.

The main result of the first paper was the following.

**Theorem 1.2** (From [MS10]) *The principal graph of any subfactor of index between 4 and 5 is a translate of one of an explicit finite list of graph pairs, which we call the* <u>vines</u>*, or is a translated extension of one of the following graph pairs, which we call the* <u>weeds</u>*.*

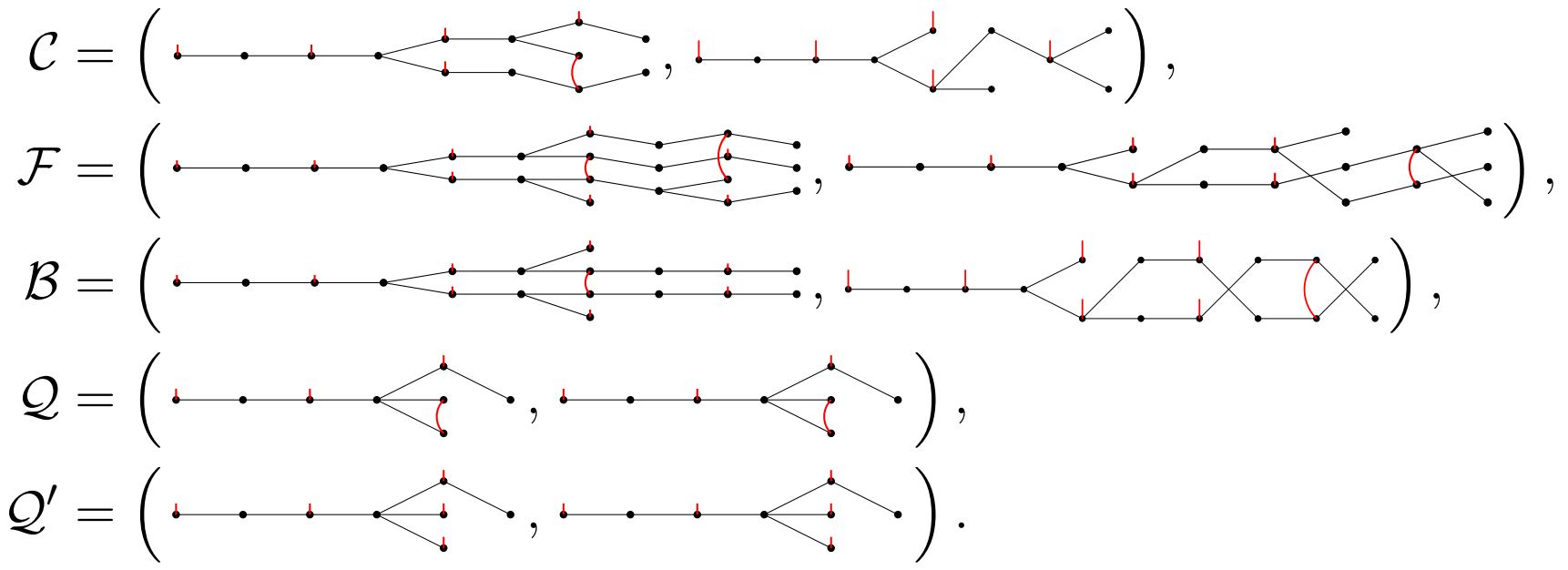

(As in [MS10], the trivial bimodule always appears as the leftmost vertex of a principal graph. Dual pairs of bimodules are indicated in red at even depths, and by matching up vertices on the two graphs at corresponding heights at odd depths.)

The main result of this paper is that three of the above weeds do not yield principal graphs of subfactors. Our technique is to use two stronger versions of the triple point obstruction. One is Jones's which comes from "quadratic tangles" planar algebra techniques [Jon03]. The other we call "the triple-single obstruction" and is proved by more traditional connections arguments. The weed $\mathcal{C}$ is ruled out by the quadratic tangles triple point obstruction, the weed $\mathcal{B}$ is ruled out by a connections argument and some combinatorial calculations, and the weed $\mathcal{F}$ is ruled out by a combination of these techniques.

**Theorem 1.3** *There are no subfactors, of any index, with principal graphs a translated extension of the pairs*

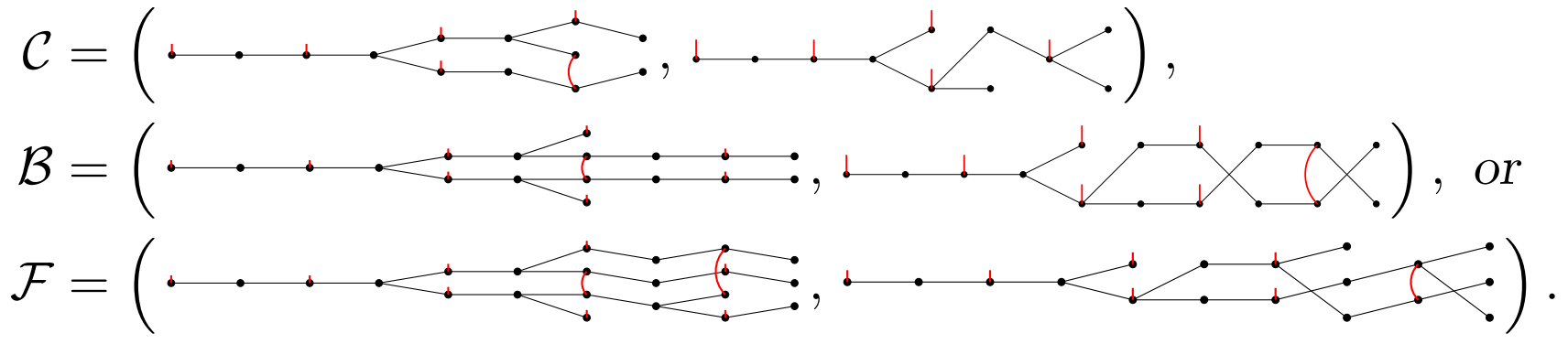

*Remark.* The triple point obstructions in this paper are also capable of ruling out a large subset of the vines described in the first paper [MS10]. We illustrate this by eliminating all but one graph in the Asaeda–Haagerup family, which is an unpublished result of Haagerup. However, applying these techniques requires a certain amount of work for each vine and do not suffice to eliminate all vines. Happily, there is a uniform arithmetic approach, which works for all vines, based on [CMS10]. A later paper in this series [PT10] uses that technique to reduce the vines to a finite set of graphs.

The structure of this paper is as follows. Section 2 quickly recalls some background about the "annular structure" of subfactors (coming from planar algebras) and connections on pairs of principal graphs (the existence of a connection is a necessary but not sufficient condition for a pair of graphs to be principal graphs of a subfactor.) Section 3 outlines three triple point obstructions, which are techniques that rule out (as principal graphs) certain graphs containing triple points. Section 4 applies these triple point obstructions to rule out the Asaeda-Haagerup vine, and the weeds $\mathcal{C}$, $\mathcal{B}$, and $\mathcal{F}$. We include an appendix with calculations of graph norms for certain infinite graphs, and closely related cousins; these are required in §4.4.

Bundled with the arXiv sources of this article are two Mathematica notebooks, `Crab.nb` and `FSM.nb`, which contain all relevant calculations for what follows. These rely on a Mathematica package called the FusionAtlas, written by the authors. See [MS10] for a terse tutorial on its use. Note that in this paper, unlike in several of the other papers in the series, every calculation can be easily checked by hand and thus this paper does not use a computer in an essential way. A typical calculation in this paper involves solving a system of a dozen or so linear equations or multiplying several polynomials in a single variable. Nonetheless we have included notebooks which perform these calculations because computer calculations are easier to check and less prone to minor errors than calculations by hand.

We would like to thank Vaughan Jones for helpful conversations and for hosting several "Planar algebra programming camps" where much of this work was done. In addition, some of the research was done at Canada/USA Mathcamp. We'd like to thank Zhengwei Liu for the argument in Lemma 3.3. During this work, Scott Morrison was at Microsoft Station Q at UC Santa Barbara and at the Miller Institute for Basic Research at UC Berkeley, David Penneys was supported by UC Berkeley's Geometry, Topology, and Operator Algebras NSF grant EMSW21-RTG, Emily Peters was in part at the University of New Hampshire and in part supported by an NSF Postdoctoral Fellowship at MIT, and Noah Snyder was supported in part by RTG grant DMS-0354321 and in part by an NSF Postdoctoral Fellowship at Columbia University. We would also like to acknowledge support from the DARPA HR0011-11-1-0001 grant.

# 2 Background

## 2.1 Annular Temperley-Lieb

The goal of this paper is to describe various triple-point obstructions, and apply these to graph pairs having annular multiplicities $*10$. We rapidly recall the language from [Jon01, Jon03] to make sense of this statement and put it in context.

A subfactor is called <u>$n$-supertransitive</u> if up to the $n$-box space its planar algebra is just Temperley-Lieb. Equivalently, a subfactor is $n$-supertransitive if and only if the principal graph up to depth $n$ is $A_{n+1}$.

Any planar algebra is a module for the annular Temperley-Lieb algebra, and as such decomposes into irreducible modules. The theory of annular Temperley-Lieb modules is laid out in Graham-Lehrer [GL98] (and in Jones [Jon01], where the idea to apply annular Temperley-Lieb theory to planar algebras appears). Each such module is cyclic, generated by a 'lowest weight vector' (that is, an irreducible submodule of a planar algebra $\mathcal{P}$ is a direct sum of subspaces of $\mathcal{P}_k$ closed under the action by annular Temperley-Lieb tangles; the weight of a vector in $\mathcal{P}_k$ is $k$, and for $n$ the lowest weight appearing in a submodule, the subspace of $\mathcal{P}_n$ is one dimensional).

Each such lowest weight vector with nonzero weight $n$ has a rotational eigenvalue which is an $n$-th root of unity. (Lowest weight vectors with weight 0 have instead a 'ring eigenvalue'.)

The annular multiplicities of a planar algebra are the sequence of multiplicities of lowest weight vectors, ignoring eigenvalues. A theorem of Jones [Jon01] shows that the annular multiplicities are actually determined entirely by the principal graph. Thus we can discuss the annular multiplicities of a graph pair regardless of whether it comes from a subfactor.

The 0th annular multiplicity of a subfactor planar algebra is always 1, corresponding to the empty diagram which generates Temperley-Lieb as an annular Temperley-Lieb module. If the planar algebra is $n$-supertransitive, then the next $n$ annular multiplicities are 0, because the vector spaces $\mathcal{P}_1$ through $\mathcal{P}_n$ are each no larger than their Temperley-Lieb subalgebra. An $n$-supertransitive subfactor of annular multiplicities $*10$ means that the first two annular multiplicities after the long string of $n$ zeroes are 1 and 0.

## 2.2 Connections

In this subsection we rapidly recall the theory of paragroups and biunitary connections on graphs developed by Ocneanu in [Ocn88]. We have taken our conventions and normalizations from [EK98]. This section is especially brief because we do not need the key notion of "flatness" for a connection, because all the obstructions in this paper are obstructions even to the existence of non-flat biunitary connections.

For simplicity, we will explain the theory for simply laced graphs. This is not at all necessary for the theory of connections, but will make our notation cleaner. Furthermore, any non-simply laced graph has index larger than 5. Nonetheless all the arguments in this paper work for arbitrary graphs. In particular, we can rule out non-simply laced principal graphs which start like any of our three weeds.

Suppose we have a bigraph pair $(\Gamma, \Gamma')$ (recall from [MS10] that this means $\Gamma$, $\Gamma'$ are bipartite graphs with dual data and specified root vertices, that $\Gamma$ and $\Gamma'$ have the same supertransitivity, and at each odd depth, there is a bijection called duality between vertices of $\Gamma$ and $\Gamma'$.) We can assemble a 4-partite Ocneanu graph from $\Gamma$ and $\Gamma'$:

$$\begin{array}{ccc} V_{00} & \xrightarrow{\ \Gamma\ } & V_{10} \\ \Gamma \Big| & & \Big| \Gamma' \\ V_{01} & \underset{\Gamma'}{\text{———}} & V_{11} \end{array}$$

(Here $V_{00}$ is the set of even vertices of $\Gamma$, $V_{11}$ is the set of even vertices of $\Gamma'$, and $V_{01}$ and $V_{10}$ are each the odd vertices of $\Gamma$, which are naturally identified with the odd vertices of $\Gamma'$.)

If we started with a subfactor $N \subset M$, the graph built from $N-N$, $N-M$, $M-M$ and $M-N$ bimodules under fusion with $X$ and $X^*$ is a 4-partite Ocneanu graph,

with

$$\begin{array}{ccc} V_{00} = \{N - N \text{ bimodules}\} & \overset{\Gamma}{\text{———}} & V_{10} = \{M - N \text{ bimodules}\} \\ \Gamma \Big| & & \Big| \Gamma' \\ V_{01} = \{N - M \text{ bimodules}\} & \underset{\Gamma'}{\text{———}} & V_{11} = \{M - M \text{ bimodules}\} \end{array}$$

**Example 2.1** *The 4-partite Ocneanu graph of the principal graphs of the Haagerup subfactor*

$\Big(\Gamma =$ 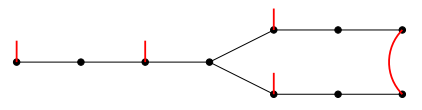 $, \Gamma' =$ 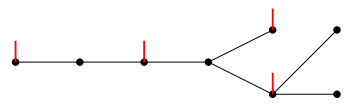 $\Big)$

*is given by*

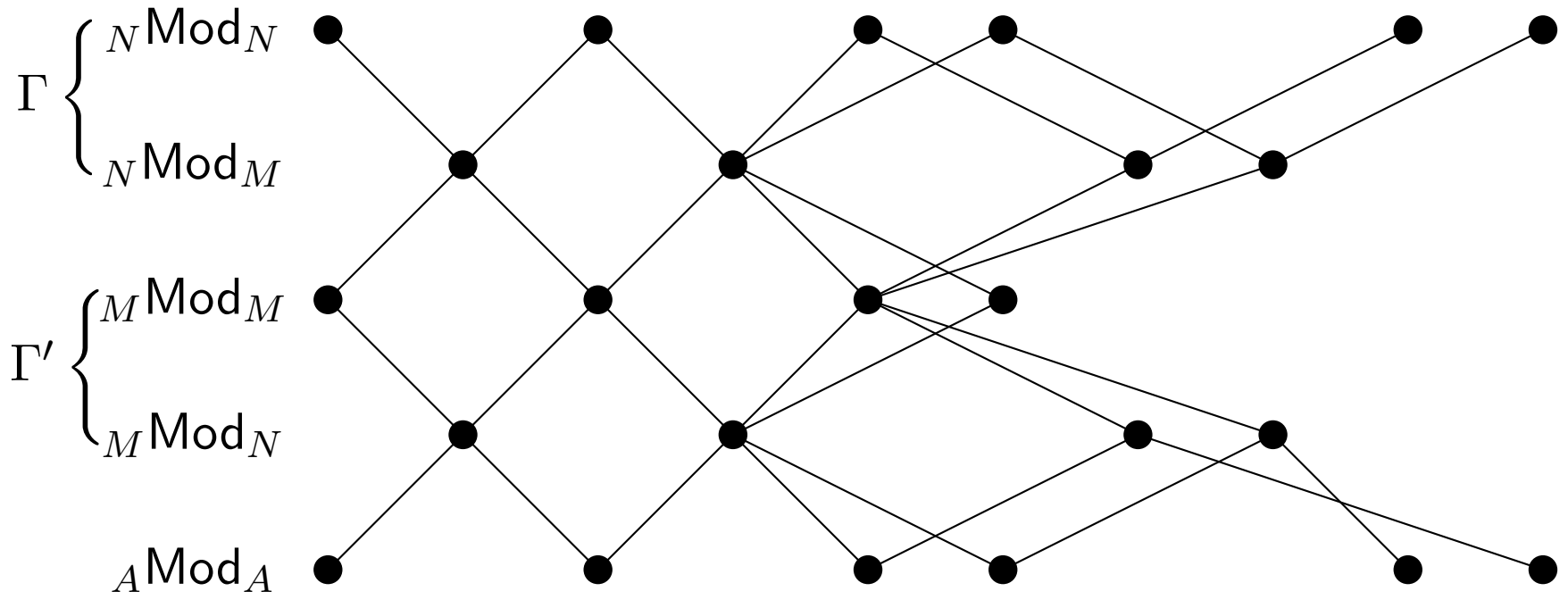

*Vertices in the first and second rows are ordered lexicographically by (depth,height) in* $\Gamma$*. Vertices in the third and fourth rows are ordered lexicographically by (depth,height) in* $\Gamma'$*.*

**Definition 2.2** *A* *<u>dimension function</u>* *on a bipartite graph* $\Gamma$ *is an eigenvector of the adjacency matrix of* $\Gamma$*, normalized so that the value of the starred vertex is* $1$*. We often think of this as an assignment of a dimension to each vertex in the graph.*

*A* *<u>dimension function</u>* *on a bigraph pair* $(\Gamma, \Gamma')$ *is pair of dimension functions, with the property that the dimensions of dual vertices are equal to each other.*

*A* *<u>dimension function</u>* *on a 4-partite graph* $\mathcal{G}$ *is an eigenvector of the adjacency matrix of* $\mathcal{G}$*, normalized so that the value of the starred vertex is* $1$*. Note that dimension functions on 4-partite graphs are in bijections with dimension functions on the corresponding bigraph pairs.*

*A dimension function is called* *<u>positive</u>* *if the dimension of each vertex is a positive real number.*

*Remark.* For finite graphs, there is a unique positive dimension function given by the Perron-Frobenius eigenvector.

**Definition 2.3** *Suppose we have a 4-partite Ocneanu graph* $\mathcal{G}$*. Then a* *<u>connection</u>* *on* $\mathcal{G}$ *is a positive dimension function* dim *and a map*

$$W : \{\textit{based loops of length 4 around } \mathcal{G}\} \to \mathbb{C}$$

*where the based loops are in one of four orders:*

$$\begin{aligned} V_{00} \to V_{10} \to V_{11} \to V_{01} \to V_{00}, \\ V_{10} \to V_{00} \to V_{01} \to V_{11} \to V_{10}, \\ V_{01} \to V_{11} \to V_{10} \to V_{00} \to V_{01}, \\ \text{or } V_{11} \to V_{01} \to V_{00} \to V_{10} \to V_{11}. \end{aligned}$$

*A connection is* <u>*biunitary*</u> *if the following properties hold.*

- <u>*Unitarity*</u>*: for all vertices $A, C$ diagonally opposite each other in $\mathcal{G}$, the matrix $W(A, -, C, -)$ is unitary; ie,*

$$\sum_D W(A,B,C,D)\overline{W(A,B',C,D)} = \delta_{B,B'}.$$

- <u>*Renormalization*</u>*: for all $A$, $B$, $C$, $D$,*

$$W(A,B,C,D) = \sqrt{\frac{\dim(B)\dim(D)}{\dim(A)\dim(C)}}\,\overline{W(B,A,D,C)}.$$

**Theorem 2.4** *A subfactor $N \subset M$ defines a biunitary connection on its 4-partite principal graph.*

**Remark 2.5** *In fact, any subfactor gives a* <u>*flat*</u> *biunitary connection. Furthermore, flatness exactly characterizes the connections which come from subfactors. On the other hand, flatness is a subtler condition which is much more difficult to check. In this paper we do not ever use the notion of flatness.*

**Corollary 2.6** *A 4-partite graph which does not have a biunitary connection is not the principal graph of a subfactor.*

**Proof of Theorem 2.4** See [EK98, ch.10-11] for a full proof; we give a brief outline in order to emphasize that this result is quite straightforward since we are not making any reference to flatness. We hope that including this sketch will make this paper more accessible to those familiar with tensor categories.

The vertices of the 4-partite principal graph are $N-N$, $N-M$, $M-M$ and $M-N$ bimodules. Let $X$ denote $M$ as an $N-M$ bimodule, ${}_N M_M$. We will write $X^{\pm}$ in tensor products to denote whichever of $X$ or $X^*$ is appropriate so that the tensor product is defined. There are $\dim(\operatorname{Hom}(Y \otimes X^{\pm}, Z))$ (or $\dim(\operatorname{Hom}(X^{\pm} \otimes Y, Z))$) edges between $Y$ and $Z$, (depending on whether the edge is horizontal or vertical). The dimension function on this graph is the bimodule dimension, namely the square root of the product of the left and right von Neumann dimensions.

For bimodules $A$ and $B$, define an inner product on $\sigma, \tau \in \operatorname{Hom}(A \otimes X^{\pm}, B)$ by $\langle \sigma, \tau \rangle = \sigma\tau^* \in \mathbb{C}$. Here we interpret $\sigma\tau^*$ as a complex number via

$$\begin{aligned} \operatorname{End}(B) &\xrightarrow{\sim} \mathbb{C} \\ \mathbf{1} &\longmapsto 1. \end{aligned}$$

Choose an orthonormal basis $\{\sigma^B_{A,X}\}$ for each space $\operatorname{Hom}_{\mathrm{N,M}}(A \otimes X, B)$, where $A$ is an irreducible $N-N$ bimodule and $B$ is an irreducible $N-M$ bimodule. Further, choose orthonormal bases $\{\sigma^B_{A,X}\}$ for the spaces $\mathrm{Hom}_{\mathrm{M,N}}(A \otimes X^*, B)$, where $A$ is an irreducible $M-M$ bimodule and $B$ is an irreducible $M-N$ bimodule. (If the graphs are simply laced, each of these spaces is one-dimensional and we are only making a choice of normalization.) Next, use these to define orthonormal bases $\{\sigma^B_{A,X}\}$ of $\mathrm{Hom}_{\mathrm{N,N}}(A \otimes X^*, B)$ or $\mathrm{Hom}_{\mathrm{M,M}}(A \otimes X, B)$ by requiring rotational invariance. That is, the Frobenius reciprocity isomorphism $\mathrm{Hom}(A \otimes X^{\pm}, B) \to \mathrm{Hom}(X^{\pm} \otimes B^*, A^*)$ should take $\sigma^B_{A,X}$ to the previously defined basis element, $\sigma^{A^*}_{X,B^*}$.

Recall that we can take adjoints of intertwiners of bimodules, giving an antilinear map $\mathrm{Hom}(X^{\pm} \otimes A, B) \to \mathrm{Hom}(B, X^{\pm} \otimes A)$. Combining this with a Frobenius reciprocity isomorphism gives the antilinear map

$$\mathrm{Hom}\left(X^{\pm} \otimes A, B\right) \to \mathrm{Hom}\left(A^* \otimes X^{\mp}, B^*\right).$$

Finally, define orthonormal bases $\{\sigma^B_{X,A}\}$ of $\mathrm{Hom}(X^{\pm} \otimes A, B)$ by requiring that $\sigma^B_{X,A}$ is sent by the above map to the previously defined basis element $\sigma^{B^*}_{A^*,X}$.

Now, $\mathrm{Hom}(X^{\pm} \otimes A \otimes X^{\pm}, C)$ has two different orthonormal bases:

$$\mathcal{L} = \{\sigma^C_{B,X}(\sigma^B_{X,A} \otimes \mathbf{1}) | B \subset X^{\pm} \otimes A \text{ and } C \subset B \otimes X^{\pm}\} \text{ and}$$
$$\mathcal{R} = \{\sigma^C_{X,D}(\mathbf{1} \otimes \sigma^D_{A,X}) | D \subset A \otimes X^{\pm} \text{ and } C \subset X^{\pm} \otimes D\}.$$

Define

$$W(A,B,C,D) = \left\langle \sigma^C_{B,X}(\sigma^B_{X,A} \otimes \mathbf{1}), \sigma^C_{X,D}(\mathbf{1} \otimes \sigma^D_{A,X}) \right\rangle$$

Then, $W(A,B,C,D)$ is the coefficient of $\sigma^C_{X,D}(\mathbf{1} \otimes \sigma^D_{A,X}) \in \mathcal{R}$ when $\sigma^C_{B,X}(\sigma^B_{X,A} \otimes \mathbf{1}) \in \mathcal{L}$ is written as a linear combination of the basis $\mathcal{R}$:

$$\sigma^C_{B,X}(\sigma^B_{A,X} \otimes \mathbf{1}) = \sum_D W(A,B,C,D) \sigma^C_{X,D}(\mathbf{1} \otimes \sigma^D_{A,X}).$$

Unitarity follows from the fact that the bases $\mathcal{L}_{A,C}$ are orthonormal; the renormalization axiom follows from the behavior of our bases $\sigma$ under duality and rotation. □

# 3 Triple Point Obstructions

## 3.1 The triple-single obstruction

Though to our knowledge the main result of this section is unpublished, we expect that some version of the following result was used by Haagerup to rule out the Asaeda-Haagerup vines beyond the first one (as we do in section 4.2).

**Theorem 3.1** (Triple-single obstruction) *Suppose we have a 4-partite graph $\mathcal{G}$, and its component graphs $\Gamma$ and $\Gamma'$ have a pair of dual triple points $\beta$ and $\beta^*$ at an odd depth, one of which is adjacent to a degree-one vertex $\gamma_2$:*

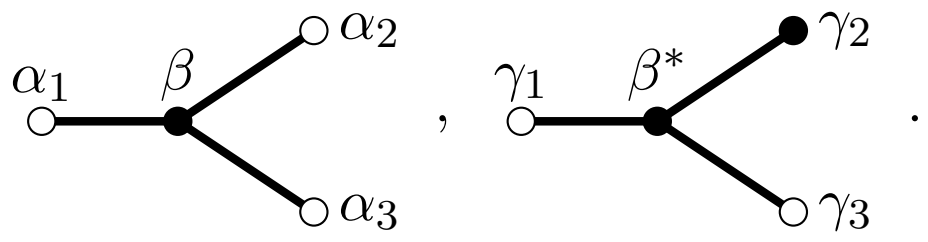

*(Observe the convention we use here: we only show a subgraph of the entire graph; if a vertex is solid, all the incident edges in the larger graph appear in the subgraph, while open vertices may be connected to other vertices in the larger graph which are not shown.)*

*Further suppose there is a biunitary connection $K$ on $\mathcal{G}$, and*

- $\dim(\alpha_1) = \dim(\gamma_1)$*;*
- *the only length-two paths (in $\mathcal{G}$) between $\alpha_1$ and $\gamma_2$ or $\gamma_3$ go through $\beta$ or $\beta^*$;*
- *the only length-two paths (in $\mathcal{G}$) between $\gamma_1$ and $\alpha_2$ or $\alpha_3$ go through $\beta$ or $\beta^*$.*

*Then*

$$|\dim(\alpha_2) - \dim(\alpha_3)| \leq K(\beta, \alpha_1, \beta^*, \gamma_1) \dim(\beta). \tag{3.1}$$

**Proof** The idea of this proof is to write down the three-by-three matrix $K(\beta, -, \beta^*, -)$; the conclusion will follow from unitarity.

Let $a_i = \sqrt{\dim(\alpha_i)}$, $b = \sqrt{\dim(\beta)}$ and $c_i = \sqrt{\dim(\gamma_i)}$. By our hypotheses, we can find the norms of all entries of $K(\beta, -, \beta^*, -)$ except three. For example, $K(\alpha_2, \beta, \gamma_2, \beta^*)$ is the sole entry of the 1-by-1 unitary matrix $K(\alpha_2, -, \gamma_2, -)$, ie a complex unit; so by the renormalization axiom, $|K(\beta, \alpha_2, \beta^*, \gamma_2)| = \frac{a_2 c_2}{b^2}$.

This gives us that, up to unspecified phases in each entry,

$$K(\beta, -, \beta^*, -) = \frac{1}{b^2} \begin{pmatrix} ? & a_1 c_2 & a_1 c_3 \\ a_2 c_1 & a_2 c_2 & ? \\ a_3 c_1 & a_3 c_2 & ? \end{pmatrix}.$$

Taking the inner products of the first two columns and dividing by $\frac{c_1 c_2}{b^2}$, (recall $c_1 = a_1$), we have

$$K(\beta, \alpha_1, \beta^*, \gamma_1) b^2 + e^{i\phi} a_2^2 + e^{i\psi} a_3^2 = 0$$

for some phases $\phi$ and $\psi$. Then by the triangle inequality, we have

$$\left| a_2^2 - a_3^2 \right| \leq K(\beta, \alpha_1, \beta^*, \gamma_1) b^2.$$

□

Although the hypotheses of this theorem seem quite stringent, they are satisfied in some interesting cases – for example, if $\beta$ is part of initial string.

**Corollary 3.2** *Suppose there is a biunitary connection on the 4-partite graph $\mathcal{G}$ with components $\Gamma$ and $\Gamma'$. Suppose $\Gamma$ and $\Gamma'$ are $(n-1)$-supertransitive (with $n$ even), there is a triple point $\beta$ at depth $n$, with dual triple point $\beta^*$, and one of the neighbors of $\beta$ or $\beta^*$ is degree-one. Then*

$$|\dim(\alpha_2) - \dim(\alpha_3)| \leq 1. \tag{3.2}$$

**Proof** The hypotheses of Theorem 3.1 are quickly verified. We find $K(\alpha_1, \beta, \gamma_1, \beta^*) = [n-1]^{-1}$ by solving for the connection along the initial $A_n$ segment (this is a quick exercise which is sketched in [EK98, p. 574-575]). This gives us that $K(\beta, \alpha_1, \beta^*, \gamma_1) = [n]^{-1}$. As $\dim(\beta) = [n]$, Theorem 3.1 implies the desired inequality. □

## 3.2 The quadratic tangles obstruction

If an $n$-supertransitive principal graph has $n$-th annular multiplicity $1$, then it begins like $D_{n+3}$ (i.e. it starts with a ‘triple point’). We define the branch factor, usually written $r$ (and $\check{r}$ for the branch factor of the dual principal graph), to be the ratio of the dimensions of the two vertices immediately past the branch point (where we take the larger divided by the smaller). If the next annular multiplicity is $0$, there are exactly three possible principal graph pairs, $(\Xi, \Xi), (\Xi, \Theta)$ or $(\Theta, \Theta)$, where

$\Xi = \cdots$

$\Theta = \cdots$

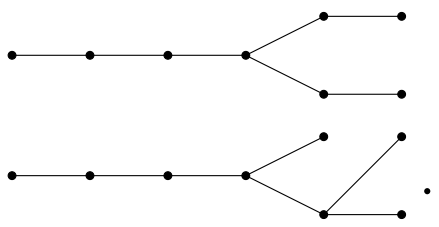

Consider now a principal graph pair with annular multiplicities $*10$, and supertransitivity $m-1$. Haagerup proved in [Haa94], using Ocneanu’s triple point obstruction (see below), that the supertransitivity must be odd, and the principal and dual principal graphs must be different. For convenience, we’ll always order the principal graph pair so the principal graph starts like the first graph above, and the dual principal graph starts like the second graph above.

An improved version of the triple point obstruction was given by Jones in [Jon03] where he also gives the following formulas for $r$, $\check{r}$ and $\lambda$, the rotational eigenvalue of the unique weight $m$ lowest weight vector.

$$r + \frac{1}{r} = \frac{\lambda + \lambda^{-1} + 2}{[m][m+2]} + 2 \tag{3.3}$$

$$\check{r} = \frac{[m+2]}{[m]} \tag{3.4}$$

The formula for $\check{r}$ follows from working out dimensions in the dual principal graph (see Example 4.2), but the formula for $r$ takes significantly more work.

Since $\lambda$ must be an $m^{\text{th}}$ root of unity, we have the following inequalities which do not involve $\lambda$:

$$-4 \leq \left(r + \frac{1}{r} - 2\right)[m][m+2] - 4 \leq 0. \tag{3.5}$$

## 3.3 Comparing the two approaches

In this section we compare the quadratic tangles triple point obstruction to the triple-single obstruction. In essence the latter is more general, while the former is stronger. In particular, we note that the triple-single obstruction is more general in the following ways.

- The triple-single obstruction has a purely local version (Theorem 3.1). That is to say, Theorem 3.1 only makes reference to the graph and the connection near the triple point. By contrast, the quadratic tangles obstruction only works for initial triple points, and thus has no local version.

- The quadratic tangles obstruction requires that the annular multiplicities be $*10$, in particular the higher depth neighbors of the 3-valent vertex on the dual graph are a 1-valent vertex and a 3-valent vertex. By contrast, Corollary 3.2 only requires that one of those vertices is 1-valent and allows the other vertex to have arbitrarily high valence.

The non-locality of the quadratic tangle obstruction is an unavoidable part of the approach: you cannot talk about rotational eigenvalues unless you know the depth. On the other hand, it's possible that the quadratic tangles approach could be modified for other initial triple-single points.

In the situation where both tests apply (namely when there's an initial triple point and annular multiplicities $*10$) the quadratic tangles test is stronger. In particular, the quadratic tangles inequality is equivalent to the triple-single inequality, but the quadratic tangles equality is stronger than both inequalities. More precisely, we have the following lemma which was pointed out to us by Zhengwei Liu.

**Lemma 3.3** *Let $\lambda$ be a root of unity, $m$ an integer, and that $a$ and $b$ are positive real numbers with $a+b=[m+1]$. Suppose that,*

$$\frac{a}{b}+\frac{b}{a}=\frac{\lambda+\lambda^{-1}+2}{[m][m+2]}+2,$$

*then*

$$(a-b)^2=\frac{[m+1]^2}{1+\frac{4[m][m+2]}{\lambda+\lambda^{-1}+2}}. \tag{3.6}$$

*(If $\lambda=-1$ the righthand side should be interpreted as $0$.) Furthermore, $|a-b|\leq 1$.*

**Proof** Equation 3.6 follows from straightforward algebra using the fact that $a+b=[m+1]$. The inequality then follows from $0\leq\lambda+\lambda^{-1}+2\leq 4$, and $[m+1]^2=1+[m][m+2]$. □

**Question 3.4** *Suppose that $\Gamma$ has an initial triple-single point and fixed supertransitivity. Can the inequality $|a-b|\leq 1$ of Corollary 3.2 be replaced by a finite set of possible values for $|a-b|$ by considering rotational eigenvalues?*

## 3.4 Ocneanu's triple point obstruction

For the sake of providing a thorough comparison of the available triple point obstructions, we briefly recall Ocneanu's obstruction, first described in Haagerup's paper [Haa94]. Note that the first paper of this series [MS10] has already made use of this obstruction. There, we stated stronger results than those described by Haagerup (but which are proved by exactly the same technique), and we merely repeat these here.

**Theorem 3.5** (Odd triple point obstruction) *Suppose we have a 4-partite graph $\mathcal{G}$ which comes from the bigraph pair $(\Gamma,\Gamma')$ of a subfactor $A\subset B$. Suppose $\Gamma$ and* $\Gamma'$ *have dual odd triple points* $\beta$ *and* $\beta^*$ *and there is some bijection* $\alpha_i \mapsto \gamma_i$ *between the neighbors of* $\beta$ *and the neighbors of* $\beta^*$*:*

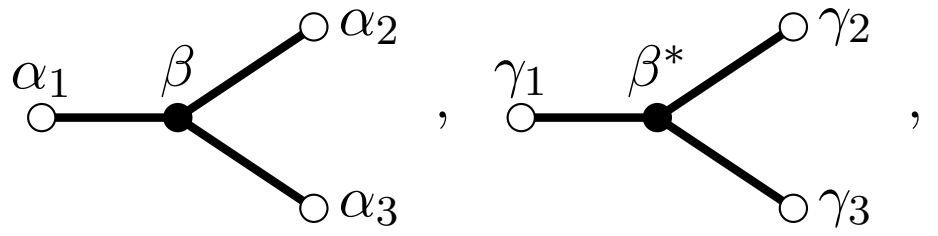

*such that*

- $\dim(\alpha_i) = \dim(\gamma_i)$ *for all* $i = 1, 2, 3$ *and*
- $\dim(\mathrm{Hom}\,(X^* \otimes \alpha_i \otimes X, \gamma_j)) = 1$ *when* $i \neq j$*.*

*Then* $[B : A] \leq 4$*.*

**Theorem 3.6** (Even triple point obstruction) *Suppose we have a* $4$*-partite graph* $\mathcal{G}$ *which comes from the bigraph pair* $(\Gamma, \Gamma')$ *of a subfactor* $A \subset B$*. Suppose* $\Gamma$ *and* $\Gamma'$ *each have a self-dual triple point at an even depth, respectively called* $\alpha$ and $\gamma$ *so that the neighbors of* $\alpha$ *are the duals of the neighbors of* $\gamma$*:*

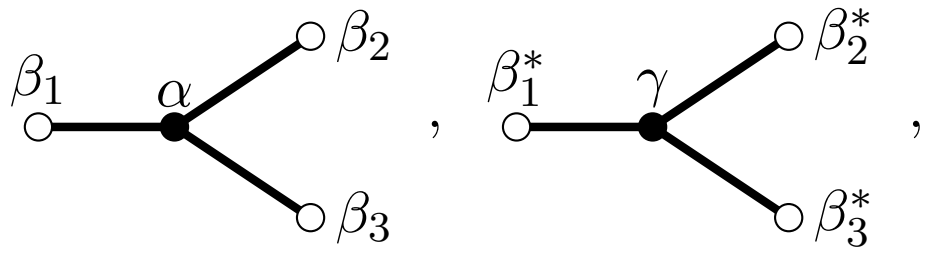

*such that*

- $\dim(\mathrm{Hom}\left(X^* \otimes \beta_i \otimes X^*, \beta_j^*\right)) = 1$ *when* $i \neq j$*.*

*Then* $[B : A] \leq 4$*.*

Typically these obstructions are used to rule out certain bigraph pairs $(\Gamma, \Gamma')$ with triple points as described above for which $\|\Gamma\|, \|\Gamma'\| > 2$. For further details, see [Haa94, MS10].

# 4 Applications

In this section we prove the main results of this paper by applying the triple point obstructions from the last section. For most of our applications we could use either the triple-single obstruction or the quadratic tangles obstructions, but typically the former will leave finitely many cases left over while using the latter approach you can apply Equation (3.3) to eliminate the exceptions.

In the first subsection we compute the dimensions of vertices in the relevant graphs. Section 4.2 then treats the Asaeda-Haagerup family, and sections 4.3, 4.4 and 4.5 eliminate the graphs $\mathcal{C}$, $\mathcal{B}$ and $\mathcal{F}$ in turn.

## 4.1 Relative dimensions of vertices

Suppose that $\Gamma$ is an $n$-translate of an extension of $\Gamma_0$ and that $\Gamma$ is the principal graph of a subfactor of index $(q+q^{-1})^2$. Because only a segment of the graph is fixed, in general the dimensions of the vertices are not determined solely by $n$ and $q$. However, in special cases we may determine these dimensions using the following three sets of equations. First, the dimension of the leftmost vertex of each graph is

$$[n+1] = \frac{q^{n+1} - q^{-n-1}}{q - q^{-1}}.$$

Second, if two vertices correspond to bimodules which are dual to each other, they must have the same dimension. Third, for each vertex $V$ which only connects to vertices which appear in the known segments of our graphs, we have an equation

$$\dim(V) = [2] \sum_{\text{edges from } W \text{ to } V} \dim(W).$$

**Remark 4.1** *Given a subfactor $A \subset B$ with principal graph $\Gamma$, not necessarily finite depth or even amenable, there are two potentially different notions of dimension for the bimodule ${}_A B_B$. One comes from the 'categorical dimension' (defined for any $A$-$B$-bimodule), which is the square root of the product of the left and right von Neumann dimensions; for ${}_A B_B$, it is $d = \sqrt{[B:A]}$. The other is the graph norm of $\Gamma$, $\|\Gamma\|$ (defined for infinite graphs in Definition A.1). If the subfactor is amenable (for example, finite-depth), then $d = \|\Gamma\|$. If the subfactor is not amenable, then $d > \|\Gamma\|$. In either case, we have that the index of the subfactor is greater than or equal to the graph norm of $\Gamma$ squared: $[B:A] \geq \|\Gamma\|^2$.*

*Even if the subfactor is non-amenable, the calculations in this section of relative dimensions by depth are still valid – the dimension relations described above allow us to express some dimensions as functions of $n$ and $q$, where $q > 1$ is such that $d = (q+q^{-1})$, even if $q$ is no longer the graph norm.*

**Example 4.2** *Suppose $(\Gamma, \Gamma')$ is an $n$-translate of an extension of*

$(\Gamma_0, \Gamma_0') = \Big($ 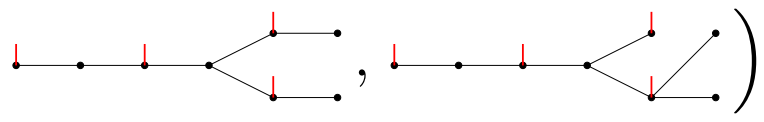 $\Big).$

*We label the vertices of $\Gamma$ and $\Gamma'$ by $V^i_{j,k}$ where $i$ is either $p$ or $d$ corresponding to either $\Gamma$ or $\Gamma'$, $j = 0, 1, \ldots, 5$ is the depth, and $k$ is the index of the vertex at that depth counting from the bottom to the top. Note that the vertices $V^i_{5,k}$ may be connected to other vertices at greater depths.*

*We cannot solve explicitly for the dimensions in terms of $n, q$ for these graphs. Rather, there is a one parameter family of solutions. We write $\alpha = \dim(V^p_{4,2})$, and then the dimensions of the vertices as functions of $n, q, \alpha$ are given by:*

$$\dim(V^p_{0,1}) = \dim(V^d_{0,1}) = \frac{q^{-n}\left(q^{2n+2}-1\right)}{q^2-1} \qquad \dim(V^p_{1,1}) = \dim(V^d_{1,1}) = \frac{q^{-n-1}\left(q^{2n+4}-1\right)}{q^2-1}$$

$$\dim(V^p_{2,1}) = \dim(V^d_{2,1}) = \frac{q^{-n-2}\left(q^{2n+6}-1\right)}{q^2-1} \qquad \dim(V^p_{3,1}) = \dim(V^d_{3,1}) = \frac{q^{-n-3}\left(q^{2n+8}-1\right)}{q^2-1}$$

$$\dim(V^p_{4,1}) = \frac{q^{-n-4}\left(\alpha q^{n+4} - \alpha q^{n+6} + q^{2n+10} - 1\right)}{q^2-1} \qquad \dim(V^p_{4,2}) = \alpha$$

$$\dim(V^d_{4,1}) = \frac{q^{-n-4}\left(q^{2n+12}-1\right)}{q^4-1} \qquad \dim(V^d_{4,2}) = \frac{q^{-n-2}\left(q^{2n+8}-1\right)}{q^4-1}$$

$$\dim(V^p_{5,1}) = \dim(V^d_{5,1}) = \frac{q^{-n-5}\left(\alpha q^{n+4} - \alpha q^{n+8} + q^{2n+12} - 1\right)}{q^2-1}$$

$$\dim(V^p_{5,2}) = \dim(V^d_{5,2}) = -\frac{q^{-n-3}\left(\alpha q^{n+2} - \alpha q^{n+6} + q^{2n+8} - 1\right)}{q^2-1}$$

**Example 4.3** *Note that in the previous example we still had a free parameter beyond $n$ and $q$, namely $\alpha$. Now we consider the case where $(\Gamma, \Gamma')$ is an $n$-translate of an extension of the graphs*

$\mathcal{C} =$ 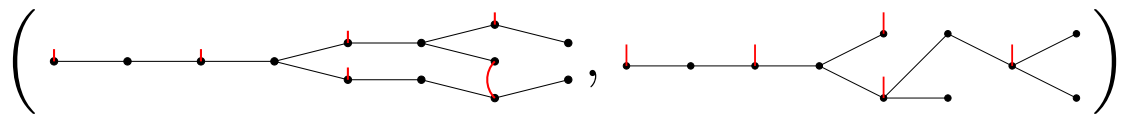

*(an extension of the previous example). Now we can solve for all the dimensions as functions of $n$ and $q$ because we have the additional equation: $[2]V^d_{5,1} = V^d_{4,1}$. The dimensions of the vertices through depth 5, as functions of $n, q$, are given by:*

$$\dim(V^p_{0,1}) = \dim(V^d_{0,1}) = \frac{q^{-n}\left(q^{2n+2}-1\right)}{q^2-1} \qquad \dim(V^p_{1,1}) = \dim(V^d_{1,1}) = \frac{q^{-n-1}\left(q^{2n+4}-1\right)}{q^2-1}$$

$$\dim(V^p_{2,1}) = \dim(V^d_{2,1}) = \frac{q^{-n-2}\left(q^{2n+6}-1\right)}{q^2-1} \qquad \dim(V^p_{3,1}) = \dim(V^d_{3,1}) = \frac{q^{-n-3}\left(q^{2n+8}-1\right)}{q^2-1}$$

$$\dim(V^p_{4,1}) = \frac{q^{-n-2}\left(q^{2n}(2q^{12}+2q^{10}+q^8) - q^4 - 2q^2 - 2\right)}{(q^2-1)(q^2+1)^3} \qquad \dim(V^p_{4,2}) = \frac{q^{-n-4}\left(q^4+q^2+1\right)\left(q^{2n+12}-1\right)}{(q^2-1)(q^2+1)^3}$$

$$\dim(V^d_{4,1}) = \frac{q^{-n-4}\left(q^{2n+12}-1\right)}{q^4-1} \qquad \dim(V^d_{4,2}) = \frac{q^{-n-2}\left(q^{2n+8}-1\right)}{q^4-1}$$

$$\dim(V^p_{5,1}) = \dim(V^d_{5,1}) = \frac{q^{-n-3}\left(q^{2n+12}-1\right)}{(q^2-1)(q^2+1)^2}$$

$$\dim(V^p_{5,2}) = \dim(V^d_{5,2}) = \frac{q^{-n-5}\left(q^{2n}(q^{16}-q^{12}-q^{10}) + q^6 + q^4 - 1\right)}{(q^2-1)(q^2+1)^2}$$

*Thus, the branch factor for this principal graph as a function of $n$ and $q$ is*

$$r(n,q) = \frac{\left(q^4+q^2+1\right)\left(q^{2n+12}-1\right)}{q^2\left(q^{2n}(2q^{12}+2q^{10}+q^8) - q^4 - 2q^2 - 2\right)}.$$

*We cannot hope for a connections argument to eliminate $\mathcal{C}$ as it is a 2-translate of a truncation of the principal graphs of the $A_3 * A_4$ Fuss-Catalan subfactor [BJ97].*

**Example 4.4** *Suppose $(\Gamma, \Gamma')$ is an $n$-translate of an extension of*

$\mathcal{F} =$ 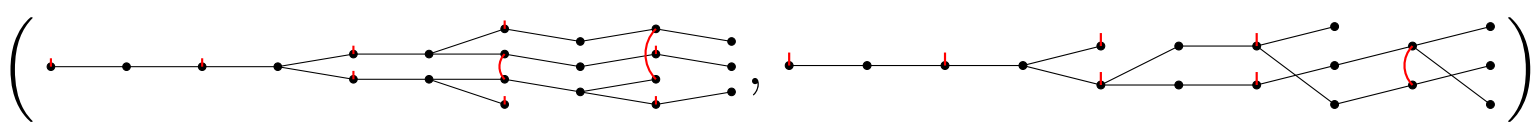 .

*Now we can solve for all the dimensions, but we need to go all the way out to depth 7 to do so. Once we do this, we find*

$$\dim(V^p_{4,1}) = \frac{q^{-4-n}\left(-1 - q^2\left(1+q^2\right)\left(2+q^2\right)\left(1+q^4\right) + q^{2(5+n)}\left(1+3q^2+3q^4+3q^6+2q^8+q^{10}\right)\right)}{\left(1+q^2\right)^3\left(-2+3q^2-3q^4+2q^6\right)}$$

$$\dim(V^p_{4,2}) = \frac{q^{-4-n}\left(-1 - q^2\left(3+2q^2+2q^4+2q^6+q^8\right) + q^{2(5+n)}\left(1+q^2\left(1+q^2\right)\left(2+2q^4+q^6\right)\right)\right)}{\left(1+q^2\right)^3\left(-2+3q^2-3q^4+2q^6\right)}.$$

*Thus, the branch factor for this principal graph as a function of $n$ and $q$ is*

$$r(n,q) = \frac{q^{2n}\left(q^{20}+3q^{18}+2q^{16}+2q^{14}+2q^{12}+q^{10}\right) - q^{10} - 2q^8 - 2q^6 - 2q^4 - 3q^2 - 1}{q^{2n}\left(q^{20}+2q^{18}+3q^{16}+3q^{14}+3q^{12}+q^{10}\right) - q^{10} - 3q^8 - 3q^6 - 3q^4 - 2q^2 - 1}.$$

*We will also need the following:*

$$\dim(V^p_{6,2}) = \frac{q^{-n-2}\left(q^4-q^3+q^2-q+1\right)\left(q^4+q^3+q^2+q+1\right)\left(q^{2n}(q^8-q^6-q^4-q^2)+q^6+q^4+q^2-1\right)}{(q-1)(q+1)\left(q^2+1\right)^3\left(2q^4-q^2+2\right)}$$

$$\dim(V^p_{5,1}) = \frac{q^{-n-1}\left(q^{2n}(q^{12}+q^6-q^2)+q^{10}-q^6-1\right)}{(q-1)(q+1)\left(q^2+1\right)^2\left(2q^4-q^2+2\right)}$$

$$\dim(V^d_{4,1}) = \frac{q^{-n}\left(q^{n+2}-1\right)\left(q^{n+2}+1\right)}{(q-1)(q+1)\left(q^2+1\right)}$$

$$\dim(V^d_{5,2}) = \frac{q^{-n-1}\left(q^{2n}(q^{12}+q^{10}-q^8-q^4-q^2)+q^{10}+q^8+q^4-q^2-1\right)}{(q-1)(q+1)\left(q^2+1\right)^2\left(2q^4-q^2+2\right)}$$

**Example 4.5** *Suppose* $(\Gamma, \Gamma')$ *is an* $n$*-translate of an extension of*

$\mathcal{B} = \Bigg($ 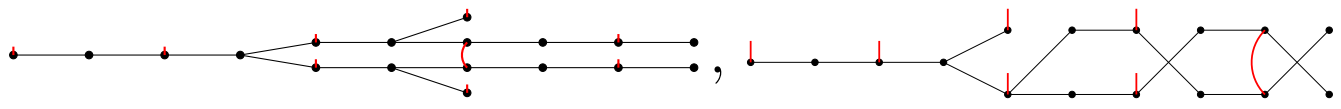 $\Bigg),$

*then we can solve for the dimensions because in the principal graph, at depth six, we have a duality between two vertices on different branches. This implies that the dimensions are the same on both branches (for those depths at which the branches remain symmetric to each other). The dimensions through depth* $6$*, as functions of* $n, q$*, are given by:*

$$\dim(V^p_{0,1}) = \dim(V^d_{0,1}) = [n+1] \qquad \dim(V^p_{1,1}) = \dim(V^d_{1,1}) = [n+2]$$

$$\dim(V^p_{2,1}) = \dim(V^d_{2,1}) = [n+3] \qquad \dim(V^p_{3,1}) = \dim(V^d_{3,1}) = [n+4]$$

$$\dim(V^p_{4,1}) = \dim(V^p_{4,2}) = \frac{[n+5]}{2} \qquad \dim(V^d_{4,1}) = \frac{[n+6]}{[2]} \qquad \dim(V^d_{4,2}) = \frac{[n+4]}{[2]}$$

$$\dim(V^p_{5,1}) = \dim(V^p_{5,2}) = \dim(V^d_{5,1}) = \dim(V^d_{5,2}) = \frac{[n+6]-[n+4]}{2}$$

$$\dim(V^p_{6,1}) = \dim(V^p_{6,4}) = \frac{[n+6]-[n+4]}{2[2]}$$

$$\dim(V^p_{6,2}) = \dim(V^p_{6,3}) = \dim(V^d_{6,1}) = \dim(V^d_{6,2}) = \frac{[n+8]-[n+6]-[n+4]-[n+2]}{2[2]}$$

*Here the branch factor* $r(n,q)$ *is equal to one, because of the duality between the branches in the principal graph. Therefore, Inequalites 3.2 and 3.5 (and indeed Equation 3.3 with* $\lambda = -1$*) always hold for translations and extensions of these graphs, and none of these triple point obstructions can eliminate this weed. Instead, this weed is eliminated in Section 4.4, based on an argument about biunitary connections.*

## 4.2 Eliminating the Asaeda-Haagerup vine

We give a proof below of an unpublished result of Haagerup stated in [Haa94].

**Theorem 4.6** *There is no biunitary connection on the* $4$*-partite graph coming from any positive translate of the Asaeda- Haagerup principal graph pair*

$\Bigg($ 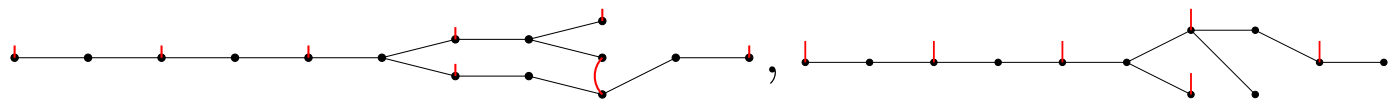 $\Bigg).$

**Proof** Suppose we translate the graphs by $j \geq 0$ so that the branch point is at depth $n = 5 + j$. Note the hypotheses of Corollary 3.2 are verified at the branch point. Labeling the vertices/ bimodules as in Example 4.2, we have $\dim(V^p_{5,1}) = [n]$, so

$$\dim(V^p_{6+k,1}) = [k+1]\dim(V^p_{6,1}) - [k][n] \text{ for } 1 \leq k \leq 4.$$

As $\dim(V^p_{10,1}) = \dim(V^p_{9,1})/[2]$, we have $\dim(V^p_{6,1}) = [5][n]/[6]$. By similar reasoning, we get the first and third equality below:

$$\frac{[3]\dim(V^p_{6,2}) - [2][n]}{2} = \dim(V^p_{8,2}) = \dim(V^p_{8,1}) = \frac{[3][5]}{[6]}[n] - [2][n],$$

and the second equality comes from duality. This means

$$\dim(V^p_{6,2}) = \left(2\frac{[5]}{[6]} - \frac{[2]}{[3]}\right)[n] = \frac{[5]+1}{[6]}[n].$$

Now by Corollary 3.2, a biunitary connection can only exist if

$$|\dim(V^p_{6,2}) - \dim(V^p_{6,1})| = \frac{[n]}{[6]} \leq 1,$$

which implies the result. □

## 4.3 Eliminating $\mathcal{C}$

In this section, we use the quadratic tangles test 3.2 to rule out principal graphs which are translated extensions of $\mathcal{C}$.

**Proposition 4.7** *Any subfactor with principal graphs a translated extension of the pair*

$\mathcal{C} =$ 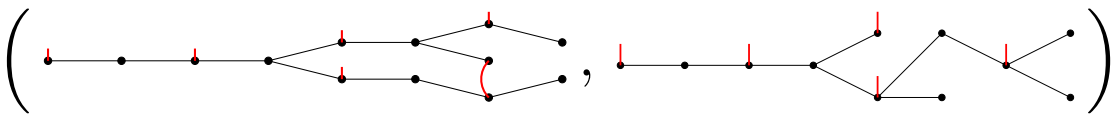

*must have index at most* $3+\sqrt{3}$.

**Proof** Suppose a subfactor exists with principal graphs an extension of the pair translated by $n \in 2\mathbb{Z}_{\geq 0}$, and let $(q+q^{-1})^2$ be the index. Plugging the branch factor

$$r(n,q) = \frac{\left(q^4+q^2+1\right)\left(q^{2n+12}-1\right)}{q^2\left(q^{2n}(2q^{12}+2q^{10}+q^8)-q^4-2q^2-2\right)}$$

calculated in Example 4.3 into Inequality (3.5) (with $m = n+4$), we get the following inequality:

$$\begin{aligned} q^{-2n-10}&\left(q^{n+5}-1\right)^2\left(q^{n+5}+1\right)^2 \times \\ &\left(q^{n+10}-q^{n+8}-q^{n+6}-q^{n+4}-q^6-q^4-q^2+1\right) \times \\ &\quad\left(q^{n+10}-q^{n+8}-q^{n+6}-q^{n+4}+q^6+q^4+q^2-1\right) \times \\ &\qquad(q-1)^{-2}(q+1)^{-2}\left(q^2-q+1\right)^{-1}\left(q^2+q+1\right)^{-1} \times \\ &\qquad\quad\left(q^{2n+8}+2q^{2n+10}+2q^{2n+12}-q^4-2q^2-2\right)^{-1} \leq 0. \end{aligned}$$

All but the two longest factors in the numerator above (namely the factors on the second and third lines) are positive for all $q > 1$. By Remark 4.1, after computing the graph norm, we see that any translated extension of the pair must satisfy $q > 1.4533$, so $q^{10} - q^8 - q^6 - q^4 > 0$, and

$$q^n \left(q^{10} - q^8 - q^6 - q^4\right) + q^6 + q^4 + q^2 - 1 \geq 0.$$

We conclude that Inequality (3.5) is satisfied if and only if

$$q^n \left(q^{10} - q^8 - q^6 - q^4\right) - q^6 - q^4 - q^2 + 1 \leq 0. \tag{4.1}$$

Note that the left hand side only increases as $n$ increases, so we examine the case $n = 0$. The largest root of

$$q^{10} - q^8 - 2q^6 - 2q^4 - q^2 + 1$$

is the positive $q$ such that $(q + q^{-1})^2 = 3 + \sqrt{3}$. Hence the index must be less than or equal to $3 + \sqrt{3}$. □

*Remark.* At this point, we could appeal to Haagerup's classification to index $3+\sqrt{3}$ to completely rule out all of these graphs. Since the published proof of his classification only covered the range up to index $3 + \sqrt{2}$, for the sake of completeness we eliminate these graphs in §4.3.

**Proposition 4.8** *Any subfactor with principal graphs a translated extension of the pair*

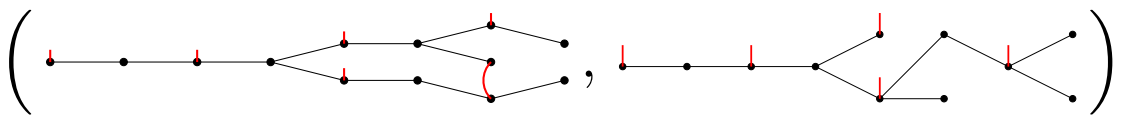

*with index less than $3 + \sqrt{3}$ is in fact a translate of one of the following graphs*

(1) 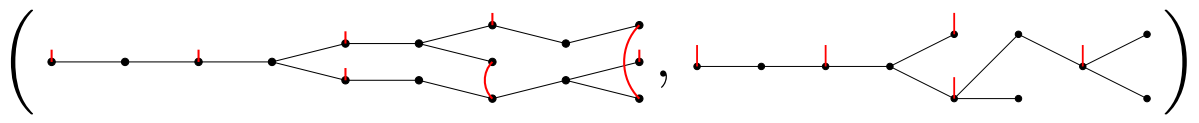

(2) 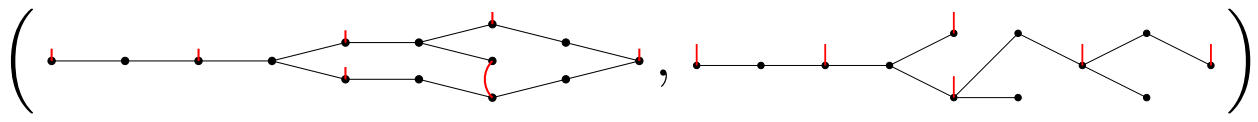

(3) 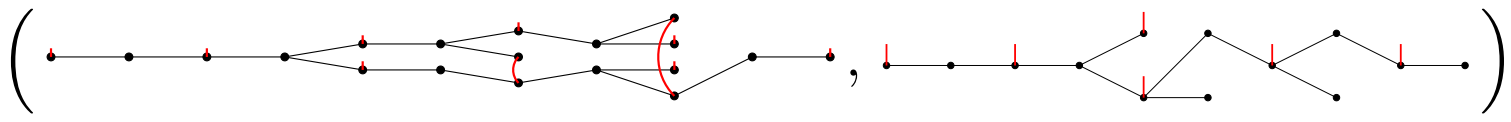

**Proof** We run the odometer, as in [MS10], and find that it terminates after two steps. (Since the index bound is low, this computation can be easily verified by hand without using a computer. In particular, the index bound forces the dual graphs to have a particularly simple form.) The four weeds considered are shown in Figure 1. Only the weed labelled 2 satisfies the associativity test, giving case 2 above. We next consider all the graphs obtained by extending one graph of a weed, staying below index $3 + \sqrt{3}$ and satisfying the associativity test. The weeds at depth $+0$ and depth $+2$ each produce exactly one such graph, giving cases 1 and 3 above.

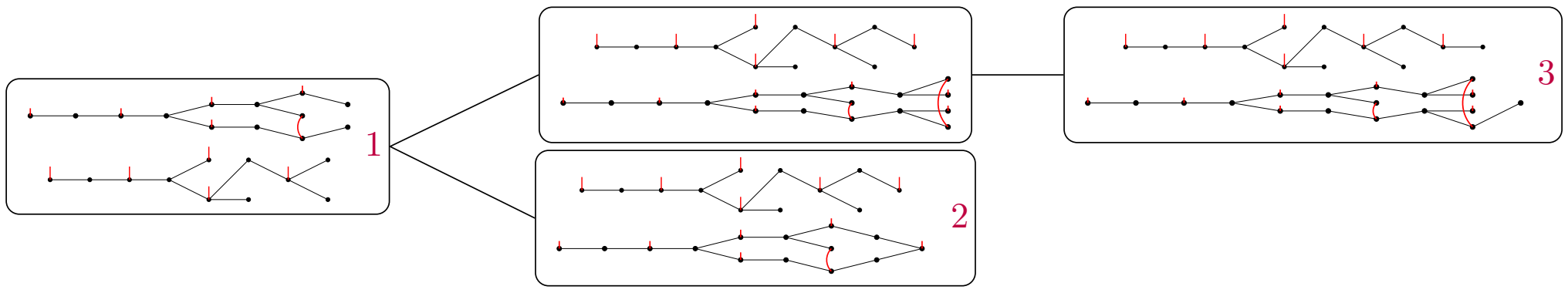

Figure 1: Running the odometer for Proposition 4.8.

□

**Proposition 4.9** *There are no subfactors with principal graphs a translation of the following pairs:*

(1) 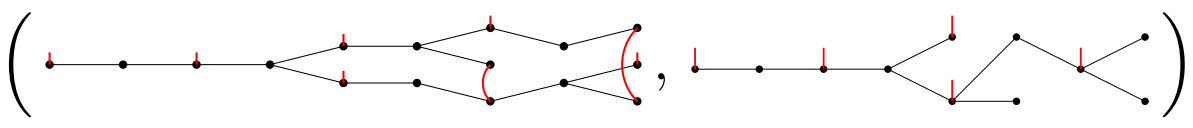

(2) 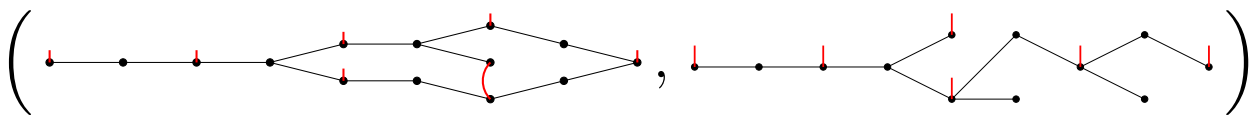

(3) 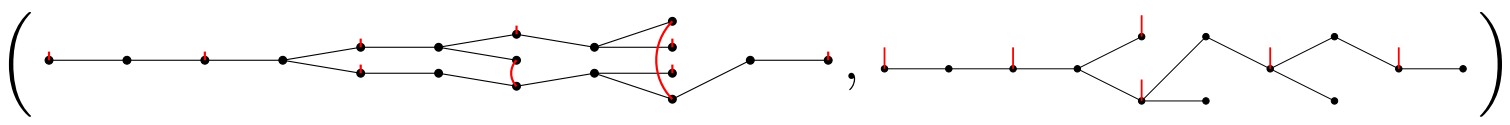

**Proof** Recall from above that for a subfactor with principal graphs a translation by $n$ of one of the above pairs and index $(q+q^{-1})^2$, we must have that $n, q$ satisfy Inequality 4.1 (which we recall for the reader's convenience):

$$q^n \left(q^{10} - q^8 - q^6 - q^4\right) - q^6 - q^4 - q^2 + 1 \leq 0.$$

For all three cases, $q > 1.4817$ by Remark 4.1, so once again

$$q^{10} - q^8 - q^6 - q^4 > 0,$$

and the left hand side of Inequality 4.1 only increases as $n$ increases. Setting $n = 2$, we have that the largest root of

$$q^{12} - q^{10} - q^8 - 2q^6 - q^4 - q^2 + 1$$

is smaller than $1.45 < 1.4817$, so this expression is always positive. Thus there cannot be subfactors with principal graphs a translation by $n \geq 2$ of any of the above pairs.

Finally, to check that these three possibilities cannot occur as principal graphs with translation $n = 0$, we note that for each case, the dimension of the lower vertex at depth 4 is not an algebraic integer. The appropriate information is contained in the table below:

| graph | minimal polynomial of dimension of vertex |
|---|---|
| 1 | $5x^3 - 16x^2 - 15x + 1$ |
| 2 | $3x^5 - 19x^4 + 25x^3 + 18x^2 - 25x - 13$ |
| 3 | $2x^2 - 6x - 9$ |

□

**Theorem 4.10** *There are no subfactors, of any index, with principal graphs a translated extension of the pair*

$\mathcal{C} =$ 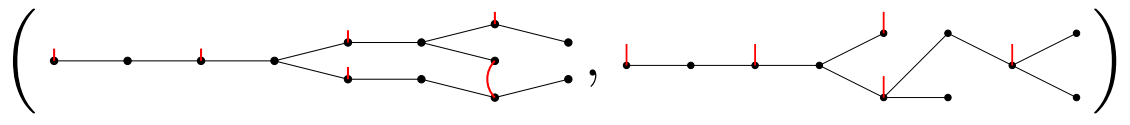 .

**Proof** The result is now an immediate consequence of Propositions 4.7, 4.8, and 4.9. □

## 4.4 Eliminating $\mathcal{B}$

In this section we show that graph pairs coming from the weed $\mathcal{B}$ cannot have a connection. First, we show that in order for a connection to exist the index must take on a particular value. Then we show that this value of the index is the square of the graph norm of the infinite depth graph $\Gamma_{n,\infty}$ (see below). Some graph theoretic arguments, which we punt to an appendix, then show that other potential principal graphs extending $\mathcal{B}$ have graph norms which are too large to satisfy the equation coming from the existence of the connection. (The graph $\Gamma_{n,\infty}$ cannot be a principal graph for other reasons.)

**Definition 4.11** *Let*

$$\mathcal{B} = (\Gamma, \Gamma') = \Big( \quad , \quad \Big)$$

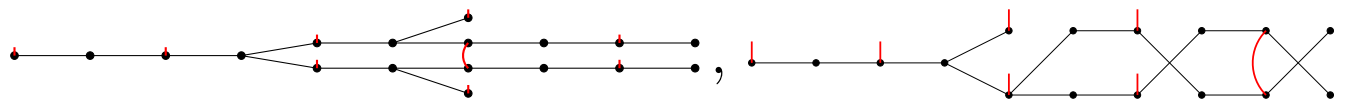

*For $n \geq 0$, let $(\Gamma_n, \Gamma'_n)$ be the translation of $(\Gamma, \Gamma')$ by $n$, so the branch point occurs at depth $n+3$. For $n \geq 0$, let $(\Gamma_{n,\infty}, \Gamma'_{n,\infty})$ be the extension of $(\Gamma_n, \Gamma'_n)$, ignoring dual data, given by adding an infinite string of vertices to each of the two maximal depth vertices. We call this the infinite simple extension of $\mathcal{B}$.*

If a principal graph pair containing $(\Gamma_n, \Gamma'_n)$ has a biunitary connection, we know what its index must be:

**Proposition 4.12** *If a subfactor has principal graph some extension of an $n$-translate of $(\Gamma, \Gamma')$, its index is $(q+q^{-1})^2$ where $q$ is the unique root greater than 1 of*

$$q^{2n+14} - 2q^{2n+12} - q^{2n+10} + q^4 + 1 = 0. \tag{4.2}$$

**Proof** The Ocneanu 4-partite graph between depths $n+3$ and $n+7$ is given by

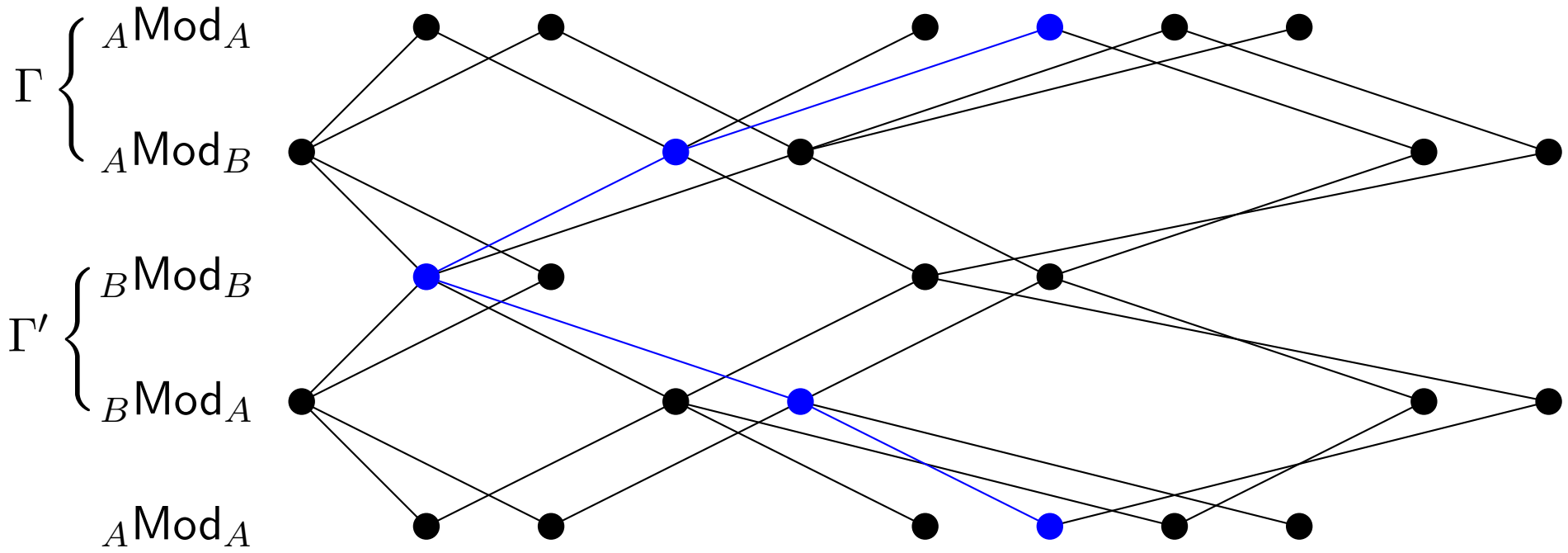

The loop $(V^p_{n+6,2}, V^p_{n+5,1}, V^d_{n+4,1}, V^d_{n+5,2})$, in blue, appears in two different 1-by-1 unitary matrices in the connection, so from the renormalization axiom we get

$$\dim(V^p_{n+6,2}) \dim(V^d_{n+4,1}) = \dim(V^p_{n+5,1}) \dim(V^d_{n+5,2}),$$

which, after expanding products of quantum numbers and canceling, gives

$$[2n+13] - [2n+11] - [2n+9] + [2n+7] - 4[n+4]^2 = 0.$$

Factoring the left hand side, we get

$$\frac{q^{-2n-10}\left(q^{2n+14}+q^{2n+10}-q^4-2q^2+1\right)\left(q^{2n+14}-2q^{2n+12}-q^{2n+10}+q^4+1\right)}{4(q-1)^2(q+1)^2\left(q^2+1\right)^2}=0.$$

Note that $q^{2n+14}+q^{2n+10}-q^4-2q^2+1>0$ for all $q>1$, so the above equality is satisfied if and only if Equation (4.2) holds. □

This norm restriction is our main tool in proving Theorem 4.13. The argument uses the fact that the restriction is satisfied by $\Gamma_{n,\infty}$; most other potential principal graphs grown from $\mathcal{B}$ have norms which are too big.

**Theorem 4.13** *There are no subfactors, of any index, with principal graphs a translated extension of the pair*

$\mathcal{B} =$ 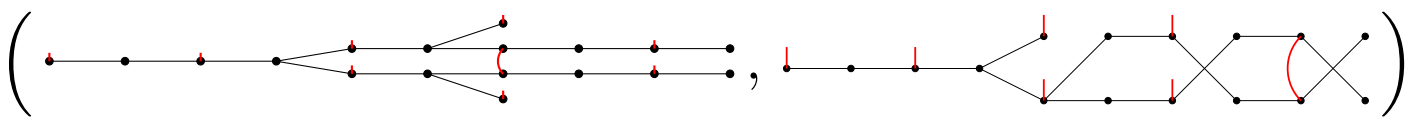

**Proof** Any translated extension of $\Gamma$, ignoring dual data, either

(A) is of the form 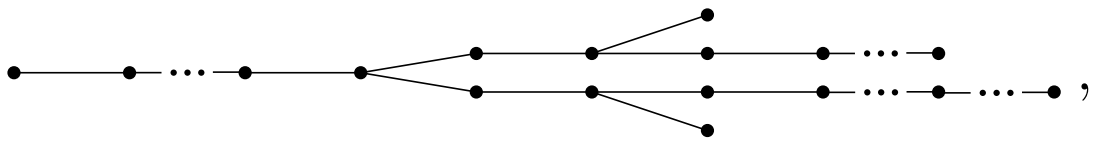

(B) is exactly $\Gamma_{n,\infty} =$ 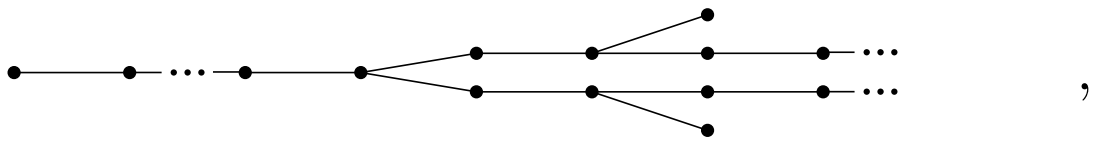

(C) contains the subgraph 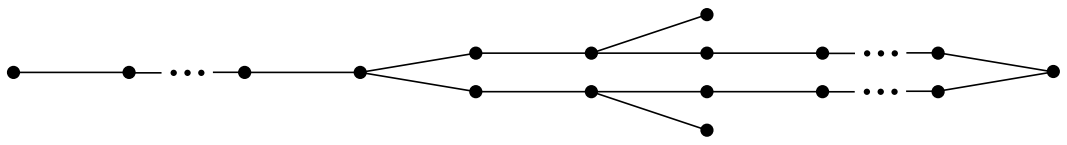 , or

(D) contains the subgraph 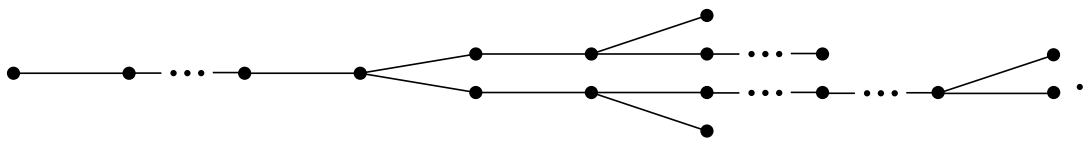

Suppose a translated extension of $\mathcal{B}$ falls in case A. Applying Lemma 4.14 inductively, we see that at depths with two vertices, the dual data continues periodically. At the depth at which either leg ends, part (1) of the lemma is contradicted.

Suppose a translated extension of $\mathcal{B}$ falls in case D. By similar reasoning, Lemma 4.14 tells us that the principal graph must actually contain

(E) 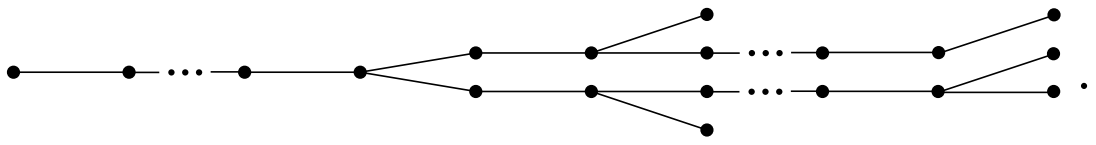

It remains to rule out the cases (B), (C) and (E). Recall that for a pair of graphs to be a principal graph pair, they must have a biunitary connection by Corollary 2.6. Proposition 4.12 tells us that if an $n$-translated extension of $\mathcal{B}$ has a biunitary connection then its norm is less than or equal to $q+q^{-1}$, where $q>1$ satisfies Equation 4.2. The graph (B) has this exactly norm, witnessed by the explicit $\ell^2$-eigenvector constructed in Lemma 4.15 below.

Using this $\ell^2$-eigenvector on (B), we can concoct a vector on either (C) or (E) giving lower bounds on the norms of these graphs. The details are given in Lemmas A.5 and A.6 in the appendix. Therefore graphs containing (C) and (E) have norms which are too large.

Finally, (B) cannot be a principal graph by Theorem 4.5 of [Pop95], or from Theorem 6.5 of [Pet09] (which applies also to infinite-depth subfactors due to [MW]). □

The proof of Theorem 4.13 relied on two lemmas, the first of which describes the allowed extensions of graphs containing a certain pattern.

**Lemma 4.14** (1) *Suppose the principal graphs* $(\Delta, \Delta')$ *of a subfactor* $A \subset B$ *between depths* $k-2$ *and* $k$ *are given by:*

$(\Delta_{odd}, \Delta'_{odd}) =$ *or*

$(\Delta_{even}, \Delta'_{even}) =$ .

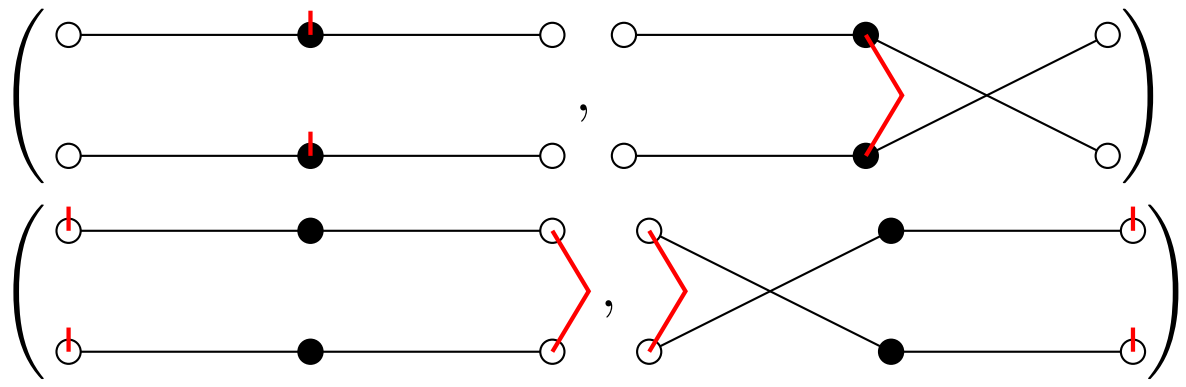

*Then each vertex at depth* $k$ *of* $(\Delta, \Delta')$ *must attach to a vertex at depth* $k+1$.

(2) *Moreover, if each vertex at depth* $k$ *of* $\Delta_{odd/even}$ *attaches to one distinct vertex at depth* $k+1$, *then each vertex at depth* $k$ *of* $\Delta'_{odd/even}$ *attaches to one distinct vertex at depth* $k+1$, *and the dual data is given by*

$(\Delta_{odd}, \Delta'_{odd}) =$ *or*

$(\Delta_{even}, \Delta'_{even}) =$

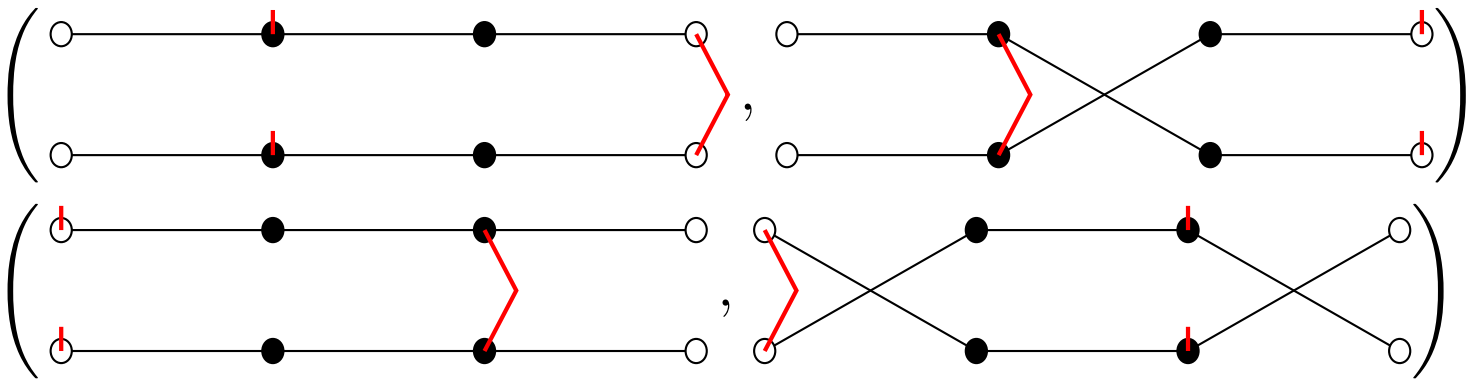

*respectively.*

**Proof** Recall that we consider bigraph pairs up to bigraph pair isomorphism and that duality of odd vertices at the same depth is given by their height in the two diagrams. We only give the proof for $(\Delta_{\text{odd}}, \Delta'_{\text{odd}})$ as the proof for $(\Delta_{\text{even}}, \Delta'_{\text{even}})$ is similar.

(1) The pair $(\Delta_{\text{odd}}, \Delta'_{\text{odd}})$ fails the associativity test unless it connects to more vertices at deeper depths. Its Ocneanu 4-partite graph, starting at depth $k-2$, is given by

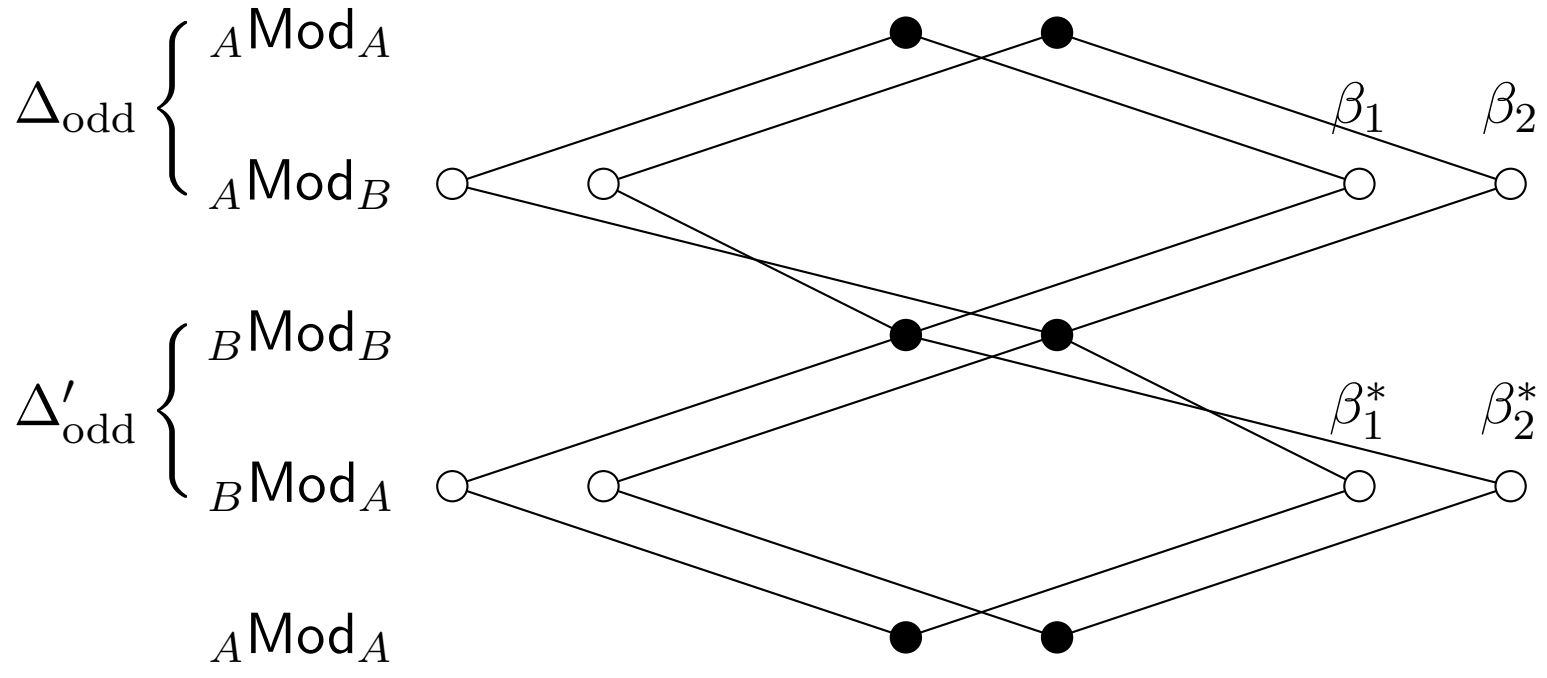

Notice there are two problems with associativity:

- There is a path from $\beta_1$ to $\beta_2^*$ through a $B-B$ bimodule, but not through an $A-A$ bimodule. The same is true for $\beta_2$ and $\beta_1^*$.
- There is a path from $\beta_1$ to $\beta_1^*$ through an $A-A$ bimodule, but not through a $B-B$ bimodule. The same is true for $\beta_2$ and $\beta_2^*$.

Hence, there must be at least one vertex at depth $k+1$ of $(\Delta, \Delta')$ attached to each vertex at depth $k$.

(2) Suppose now that each vertex $\beta_i$ at depth $k$ of $\Delta$ attaches to one distinct vertex $\alpha_i$ at depth $k+1$. To fix the first problem with associativity, the vertices at depth $k+1$ must be dual to each other:

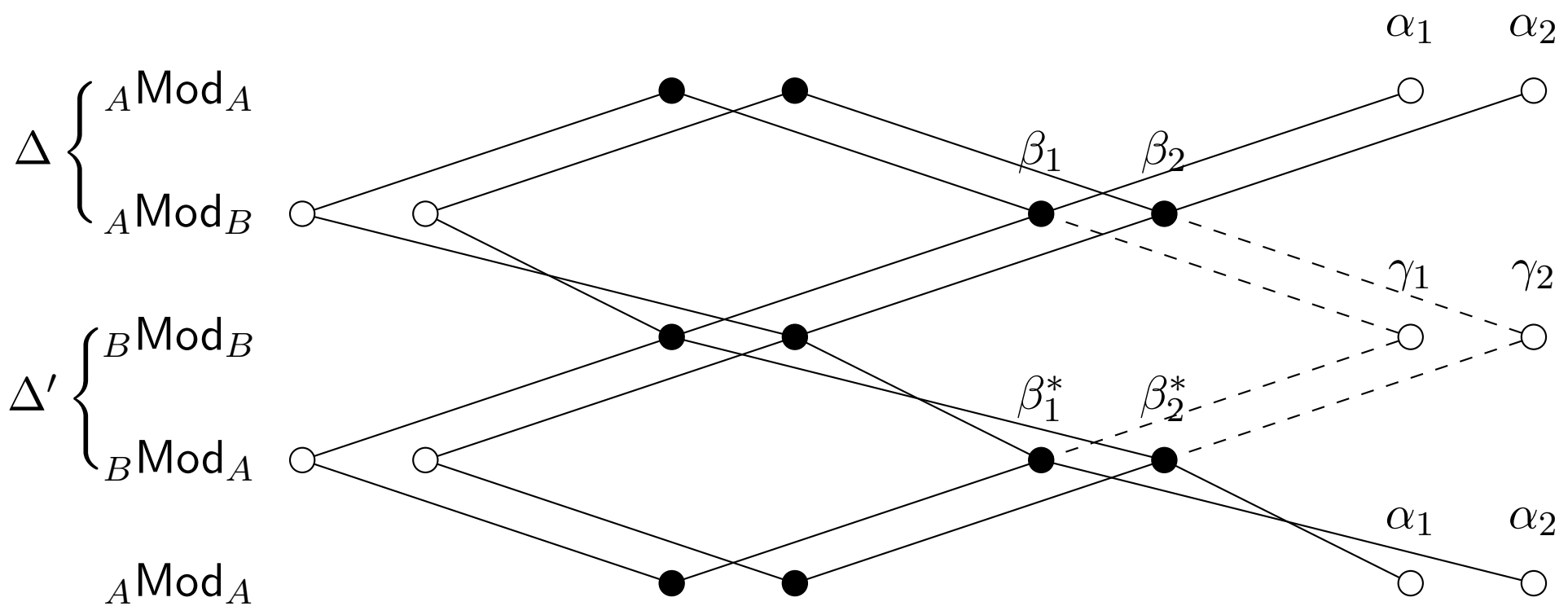

By inspection, we see that each vertex $\beta_i^*$ at depth $k$ in $\Delta'$ must be attached to a distinct, self-dual vertex $\gamma_i$ at depth $k+1$. This is the only arrangement of depth $k+1$ vertices consistent with our assumptions about $\Delta$.

□

We also needed the following lemma about the norm of $\Gamma_{n,\infty}$ (which in turn depends on some facts about norms of infinite graphs, discussed in Appendix A).

**Lemma 4.15** *Let $q$ be the unique root greater than 1 of Equation (4.2). The graph $\Gamma_{n,\infty}$ has norm $d = q + q^{-1}$, because it has a totally positive $\ell^2$-eigenvector $\mathbf{v}$ with eigenvalue $d = q + q^{-1}$.*

**Proof** We show that $\Gamma_{n,\infty}$ has a totally positive $\ell^2$-eigenvector $\mathbf{v}$, with eigenvalue $d = q + q^{-1}$. Then by Theorem A.3, $\|\Gamma_{n,\infty}\| = d$.

We will define a vector $v \in \ell^2(\Gamma_{n,\infty})$, and denote by $v_{i,j}$ its value at the $j$-th vertex at depth $i$. Recall that the branch point is at depth $n+3$. The early entries $v_{i,j}$ are the dimensions calculated in Example 4.5:

$$v_{k,1} = [k+1] \text{ for all } k \leq n+3 \qquad v_{n+4,1} = v_{n+4,2} = \frac{[n+5]}{2}$$
$$v_{n+5,1} = v_{n+5,2} = \frac{[n+6]-[n+4]}{2} \qquad v_{n+6,1} = v_{n+6,4} = \frac{[n+6]-[n+4]}{2[2]};$$

Now, Lemma A.4 tells us that if $v$ is to to be an $\ell^2$-eigenvector, the rest of the $v_{i,j}$ must be a decreasing geometric series with $q^{-1}$:

$$v_{n+6,2} = v_{n+6,3} = \frac{[n+6]-[n+4]}{2q^{-1}}$$
$$v_{n+k,1} = v_{n+k,2} = \frac{[n+6]-[n+4]}{2q^{-k+5}} \text{ for all } k \geq 7.$$

By inspection, this is an eigenvector if and only if

$$v_{n+4,1} + v_{n+6,1} + v_{n+6,2} = [2]v_{n+5,2}.$$

Expanding in terms of $q$ and factoring, we see this equation is equivalent to

$$\frac{q^{-n-4}\left(q^{2n+14} - 2q^{2n+12} - q^{2n+10} + q^4 + 1\right)}{2(q-1)(q+1)\left(q^2+1\right)} = 0,$$

□

### 4.5 Eliminating $\mathcal{F}$

In this section, we use the quadratic tangles test from §3.2 and the existence of connections on subfactor principal graphs to rule out principal graphs which are translated extensions of $\mathcal{F}$.

**Proposition 4.16** *Any subfactor with principal graphs a translated extension of the pair*

$\mathcal{F} =$ 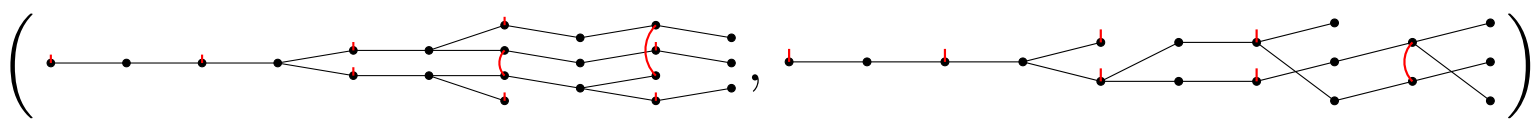

*must either*

(1) *have principal graphs translated by* $0$ *and have rotational eigenvalue* $\lambda$ *and index* $(q+q^{-1})^2$ *where* $\lambda$ *and* $q$ *are either:*

| $q$ | *minimal polynomial for* $q$ | $\lambda$ |
|---|---|---|
| 2.0106... | $x^{12} - 3x^{10} - 3x^8 - 4x^6 - 3x^4 - 3x^2 + 1$ | 1 |
| 1.8449... | $x^{36} + x^{34} - 2x^{32} - 17x^{30} - 46x^{28} - 91x^{26} - 144x^{24}$ $-197x^{22} - 233x^{20} - 246x^{18} - 233x^{16} - 197x^{14}$ $-144x^{12} - 91x^{10} - 46x^8 - 17x^6 - 2x^4 + x^2 + 1$ | $\pm i$ |

*or,*

(2) *have principal graphs translated by* $2$ *and have rotational eigenvalue* $\lambda$ *and index* $(q+q^{-1})^2$ *where* $\lambda$ *and* $q$ *are either:*

| $q$ | *minimal polynomial for* $q$ | $\lambda$ |
|---|---|---|
| 1.6341... | $x^{16} - x^{14} - 2x^{12} - 5x^{10} - 2x^8 - 5x^6 - 2x^4 - x^2 + 1$ | 1 |
| 1.6069... | $x^{52} - x^{48} - 4x^{46} - 6x^{44} - 19x^{42} - 38x^{40} - 67x^{38}$ $-98x^{36} - 139x^{34} - 178x^{32} - 218x^{30} - 238x^{28}$ $-246x^{26} - 238x^{24} - 218x^{22} - 178x^{20} - 139x^{18}$ $-98x^{16} - 67x^{14} - 38x^{12} - 19x^{10} - 6x^8 - 4x^6 - x^4 + 1$ | $\exp(\pm\pi i/3)$ |

*Remark.* Note that both cases in (1) and the first case of (2) have index strictly greater than 5.

**Proof** First note that the $q$ from any translated extension of this pair must be at least 1.5932 by Remark 4.1. Proceeding as in Proposition 4.7, the branch factor as a function of $n$ and $q$ is given by

$$r(n,q) = \frac{q^{2n}\left(q^{20} + 3q^{18} + 2q^{16} + 2q^{14} + 2q^{12} + q^{10}\right) - q^{10} - 2q^8 - 2q^6 - 2q^4 - 3q^2 - 1}{q^{2n}\left(q^{20} + 2q^{18} + 3q^{16} + 3q^{14} + 3q^{12} + q^{10}\right) - q^{10} - 3q^8 - 3q^6 - 3q^4 - 2q^2 - 1}.$$

Plugging in $r(n,q)$ to Equation (3.5), we get the following inequality:

$$q^{-2n-4}\left(q^{n+5}-1\right)^2\left(q^{n+5}+1\right)^2(q-1)^{-2}(q+1)^{-2}\times$$
$$\left(q^{2n}\left(q^{16}-q^{14}-q^{12}-q^{10}\right)+q^n\left(-2q^{14}-3q^{12}+3q^4+2q^2\right)+q^6+q^4+q^2-1\right)\times$$
$$\left(q^{2n}\left(q^{16}-q^{14}-q^{12}-q^{10}\right)+q^n\left(2q^{14}+3q^{12}-3q^4-2q^2\right)+q^6+q^4+q^2-1\right)\times$$
$$\left(q^{2n}\left(q^{20}+2q^{18}+3q^{16}+3q^{14}+3q^{12}+q^{10}\right)-q^{10}-3q^8-3q^6-3q^4-2q^2-1\right)^{-1}\times$$
$$\left(q^{2n}\left(q^{20}+3q^{18}+2q^{16}+2q^{14}+2q^{12}+q^{10}\right)-q^{10}-2q^8-2q^6-2q^4-3q^2-1\right)^{-1}\leq 0.$$

By similar analysis as above, this inequality is satisfied if and only if

$$q^{2n}\left(q^{16}-q^{14}-q^{12}-q^{10}\right)+q^n\left(-2q^{14}-3q^{12}+3q^4+2q^2\right)+q^6+q^4+q^2-1\leq 0.$$

Let $p(n,q)$ denote the left hand side. If $n\geq 4$ and $q>1$, then

$$\begin{aligned} p(n,q) &\geq q^{2n}\left(q^{16}-q^{14}-q^{12}-q^{10}\right)+q^n\left(-2q^{14}-3q^{12}\right) \\ &\geq q^{2n}\left(-2q^{10}-3q^8+q^{16}-q^{14}-q^{12}-q^{10}\right) \\ &= q^{2n+8}\left(q^8-q^6-q^4-3q^2-3\right). \end{aligned}$$

The largest root of

$$q^8-q^6-q^4-3q^2-3$$

is less than $1.5082 < 1.5932$, so there can be no subfactors with an $n$-translated extension of this pair of principal graphs for $n\geq 4$.

Now suppose we have a subfactor with principal graphs an extension of this pair of principal graphs (with no translation). Then $\lambda\in\{\pm 1,\pm i\}$ and $\lambda+\lambda^{-1}\in\{-2,0,2\}$. Solving Equation (3.3) for $q$ when $\lambda=-1$ shows that $q$ must be approximately 1.3123..., with minimal polynomial $x^8-x^6-x^4-x^2+1$. This $q$ is smaller than 1.5932 so we can ignore this case. Solving Equation (3.3) for $q$ when $\lambda\in\{1,\pm i\}$ gives the first table in the statement.

Finally, suppose we have a subfactor with principal graphs a 2-translated extension of this pair of principal graphs. Then $\lambda\in\{\pm 1,\exp(\pm 2\pi i/3),\exp(\pm\pi i/3)\}$ and $\lambda+\lambda^{-1}\in\{-2,-1,1,2\}$. Solving Equation (3.3) for $q$ when $\lambda\in\{-1,\exp(\pm 2\pi i/3)\}$ gives the cases

| $q$ | minimal polynomial for $q$ | $\lambda$ |
|---|---|---|
| 1.3453... | $x^{16}-x^{14}-2x^{10}-2x^6-x^2+1$ | $-1$ |
| 1.5203... | $x^{52}-x^{48}-4x^{46}-4x^{44}-9x^{42}-14x^{40}-21x^{38}$ $-24x^{36}-29x^{34}-36x^{32}-42x^{30}-44x^{28}-42x^{26}$ $-44x^{24}-42x^{22}-36x^{20}-29x^{18}-24x^{16}-21x^{14}$ $-14x^{12}-9x^{10}-4x^8-4x^6-x^4+1$ | $\exp(\pm 2\pi i/3)$ |

which we ignore as $q$ is too small. Solving Equation (3.3) for $q$ when $\lambda\in\{1,\exp(\pm\pi i/3)\}$ gives the second table above. □

**Proposition 4.17** *In order for a connection to exist on any extension of an $n$-translate of $\mathcal{F}$, it must have an eigenvalue $d=q+q^{-1}$ where $q$ is the unique root*

*greater than* 1 *of*

$$q^{4n}\left(q^{40}-2q^{38}+q^{36}-4q^{34}-4q^{32}-5q^{30}-4q^{28}-6q^{26}-2q^{24}-3q^{22}-q^{20}\right)$$
$$+q^{2n}\left(-q^{28}+9q^{26}+4q^{24}+11q^{22}+12q^{20}+11q^{18}+4q^{16}+9q^{14}-q^{12}\right)$$
$$-q^{20}-3q^{18}-2q^{16}-6q^{14}-4q^{12}-5q^{10}-4q^{8}-4q^{6}+q^{4}-2q^{2}+1 \quad (4.3)$$

**Proof** The proof is similar to Proposition 4.12. The Ocneanu 4-partite graph between depths $n+3$ and $n+7$ is given by

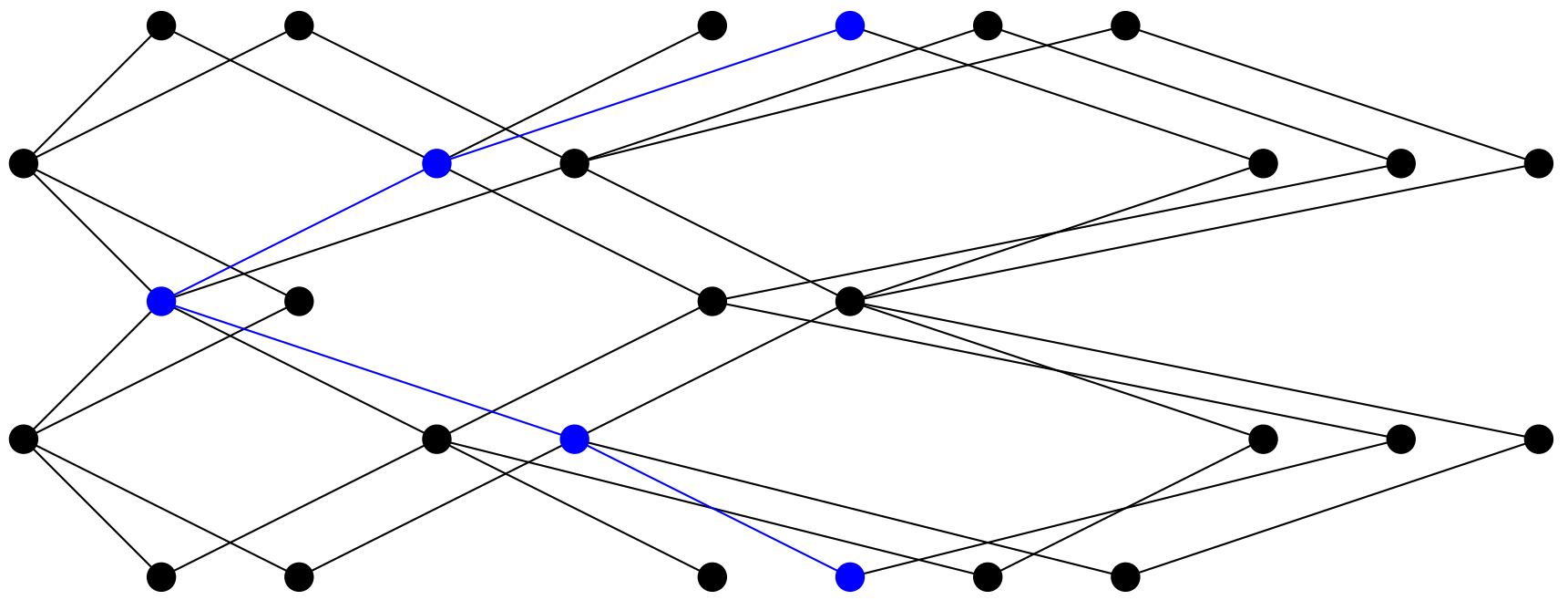

The loop $(V^p_{n+6,2}, V^p_{n+5,1}, V^d_{n+4,1}, V^d_{n+5,2})$, in blue, appears in two different 1-by- 1 unitary matrices in the connection, so the renormalization axiom gives us

$$\dim(V^p_{n+6,2})\dim(V^d_{n+4,1}) - \dim(V^p_{n+5,1})\dim(V^d_{n+5,2}) = 0.$$

These dimensions are given in Example 4.4. Equation (4.3) is obtained from the numerator after substituting and simplifying. □

**Theorem 4.18** *There are no subfactors, of any index, with principal graphs a translated extension of the pair*

$\mathcal{F} =$

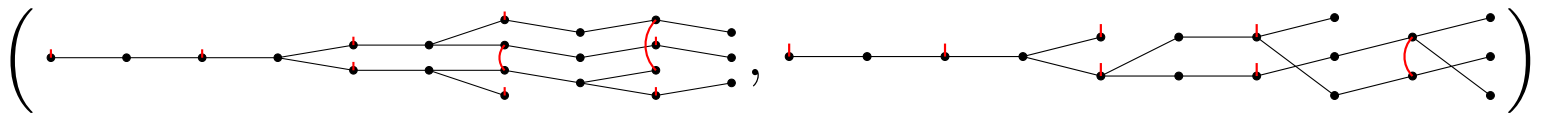

**Proof** We will show the results of Propositions 4.16 and 4.17 are mutually exclusive. Suppose we have a subfactor with principal graphs given by an $n$-translated extension of $\mathcal{F}$. By Proposition 4.16, we know that $n \in \{0, 2\}$.

Substituting $n = 0$ in Equation (4.3) and factoring, we see that $q = 1.6068...$ with minimal polynomial

$$x^{28}-x^{26}+x^{24}-6x^{22}-6x^{20}-19x^{18}-19x^{16}-27x^{14}-19x^{12}-19x^{10}-6x^{8}-6x^{6}+x^{4}-x^{2}+1,$$

a contradiction to Proposition 4.16.

Substituting $n = 2$ in Equation (4.3) and factoring, we see that $q = 1.6118...$ with minimal polynomial

$$x^{36}-2x^{34}+2x^{32}-6x^{30}-x^{28}-13x^{26}-4x^{24}-23x^{22}-10x^{20}$$
$$-24x^{18}-10x^{16}-23x^{14}-4x^{12}-13x^{10}-x^{8}-6x^{6}+2x^{4}-2x^{2}+1,$$

a contradiction to Proposition 4.16. □

## A Facts about norms for some infinite graphs

In the following $G$ will always denote a locally finite graph and $A(G)$ its adjacency matrix.

**Definition A.1** *For any graph $G$ (not necessarily finite or even locally finite), its graph norm $\|G\|$ is the operator norm of its adjacency matrix*

$$\|A(G)\| = \sup_{v \in \ell^2(\Gamma)} \frac{\|Av\|}{\|v\|}$$

*where $\ell^2(\Gamma)$ is $\ell^2$-functions on the vertices of $\Gamma$.*

**Theorem A.2** *[MW89, 4.13] Suppose subgraphs $G_n$ converge to $G$. Then $\|G_n\| \nearrow \|G\|$.*

**Theorem A.3** *If an infinite graph $G$ has an $\ell^2$-eigenvector $\mathbf{v}$ with strictly positive entries corresponding to eigenvalue $d$, then $\|G\| = d$.*

**Proof** This follows from Theorems 4.4 and 6.2 of [MW89]. See also the remark of page 183 of [Pop94]. □

**Lemma A.4** *Suppose $G$ is of the form*

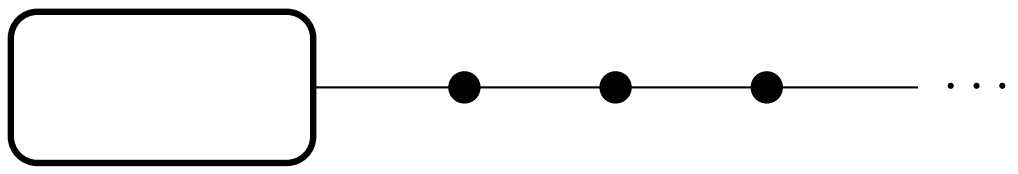

*(Here the empty rectangle indicates some arbitrary graph; outside the rectangle we have an infinite chain of edges.)*

*If $\mathbf{v}$ is a strictly positive $\ell^2$-eigenvector for $G$ with eigenvalue $d = (q + q^{-1}) > 2$, with entries*

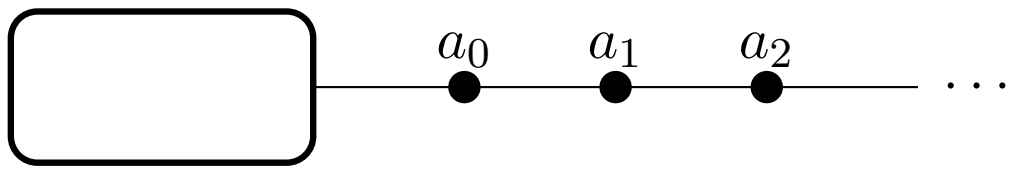

*then $a_n = q^{-n} a_0$.*

**Proof** From $a_0$ and $a_1$ and the relation $[2]a_k = a_{k-1} + a_{k+1}$, one shows inductively that $a_n = [n]a_1 - [n-1]a_0$.

For $\mathbf{v}$ to be an $\ell^2$-eigenvector, we need $a_n \to 0$. There's some $\epsilon$ such that $a_0 = (q + \epsilon)a_1$; expanding out $a_n$ using this relation, we find

$$a_n = [n]a_1 - [n-1]a_0 = a_2 \frac{-\epsilon q^{n-2} + q^{-n+3} + \epsilon q^{-n+2} - q^{-n+1}}{q - q^{-1}}.$$

Since $q > 1$ and $a_n \to 0$, we must have $\epsilon = 0$.

So we know $a_0 = qa_1$. Now $a_n = [n]a_1 - [n-1]a_0 = (q^{-1}[n] - [n-1])a_0 = q^{-n}a_0$. □

**Lemma A.5** *If*

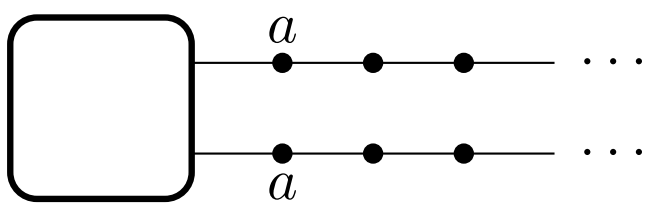

*has an $\ell^2$-eigenvector with eigenvalue $d > 2$ and the components at the two vertices marked $a$ above are each equal to $a$, then*

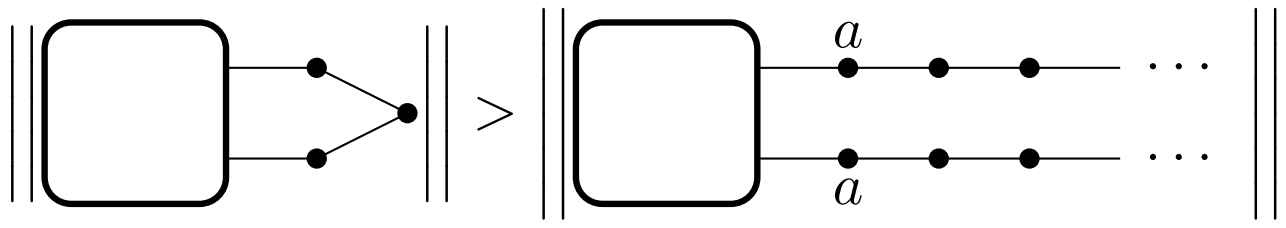

**Proof** By Lemma A.4, the $\ell^2$-eigenvector is of the form

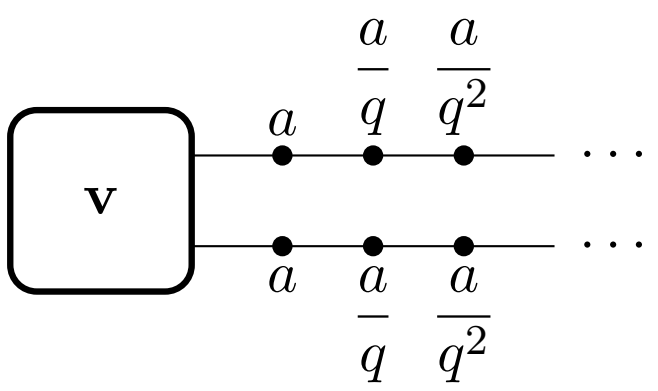

with eigenvalue $d = q + q^{-1}$.

Now, consider this vector "restricted" to the finite graph:

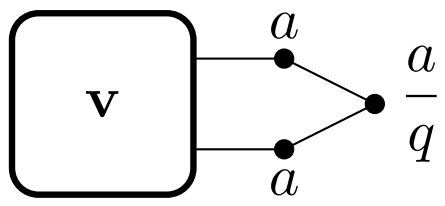

By Frobenius-Perron eigentheory, we know

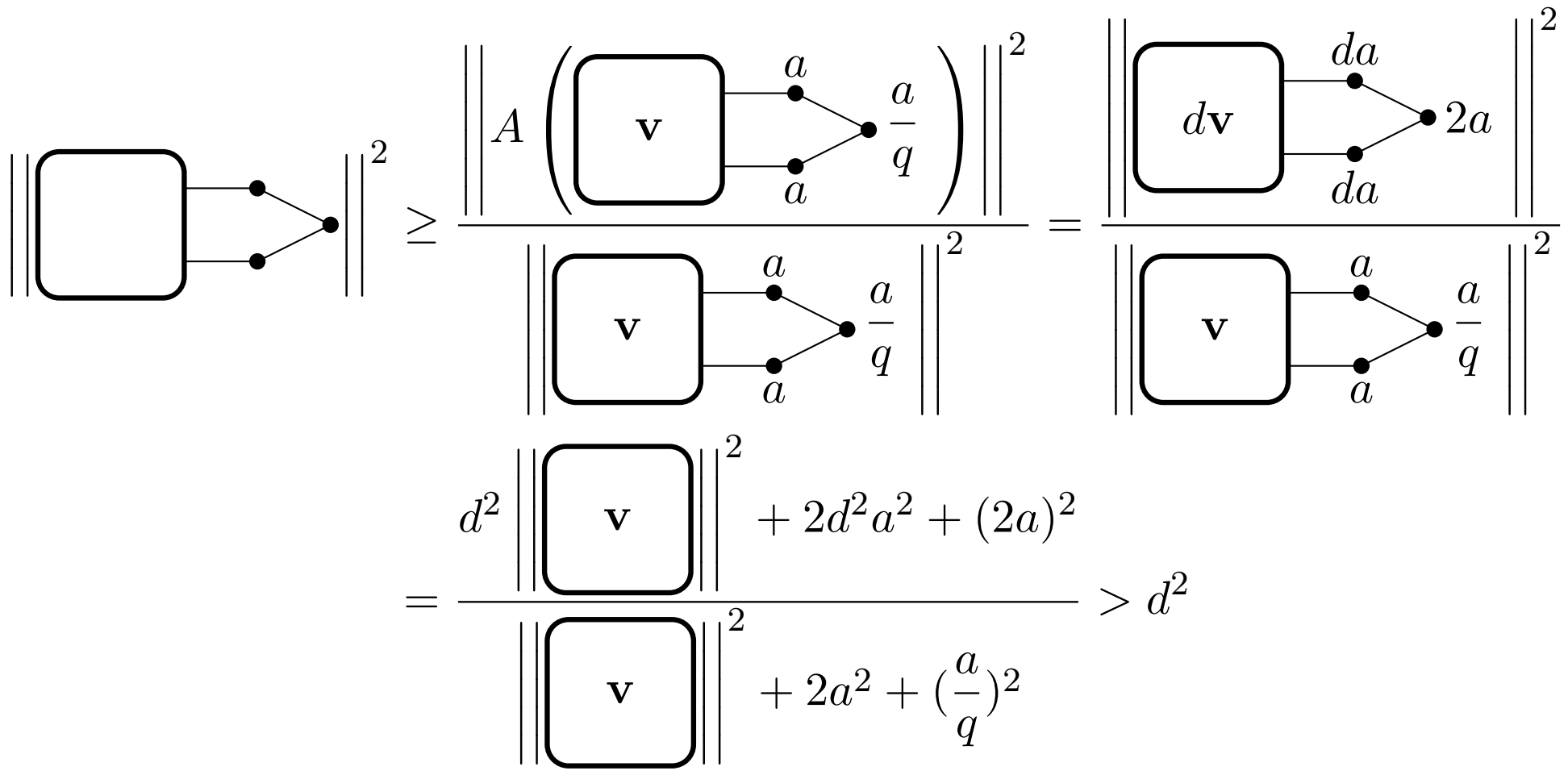

because $2 > dq^{-1}$ (since $2q > d = q + q^{-1}$). By Theorem A.3, we are finished. □

**Lemma A.6** *If*

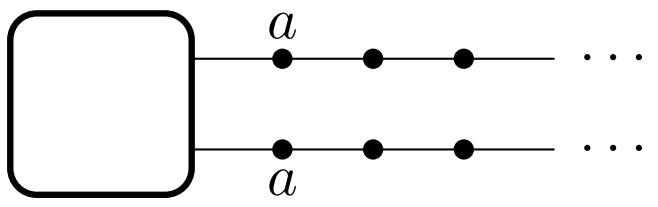

*has a strictly positive $\ell^2$-eigenvector with eigenvalue $d = q + q^{-1}$ where $2q^2 - 3 - 3q^{-2} > 0$ (which is true for $q > 1.48$), then*

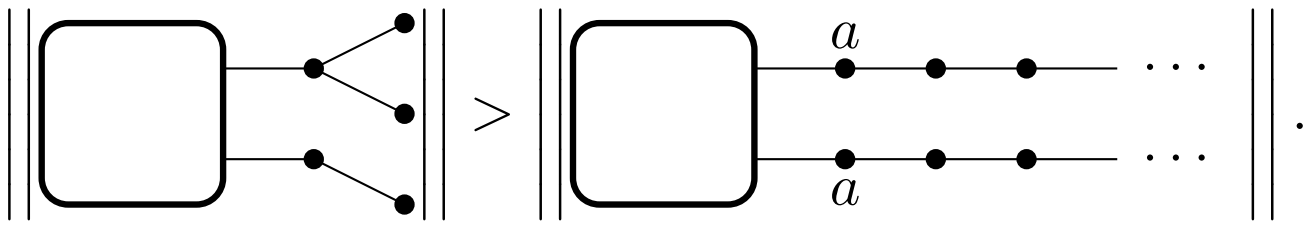

**Proof** This is similar to the proof of Lemma A.5. We consider the eigenvector "restricted" to the finite graph:

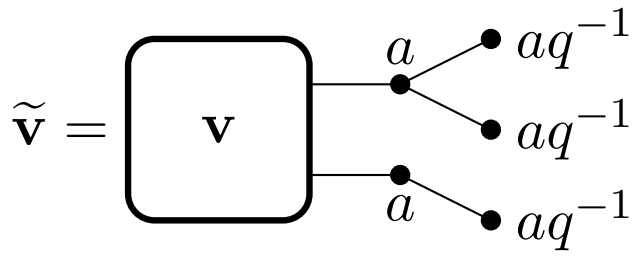

Then we have

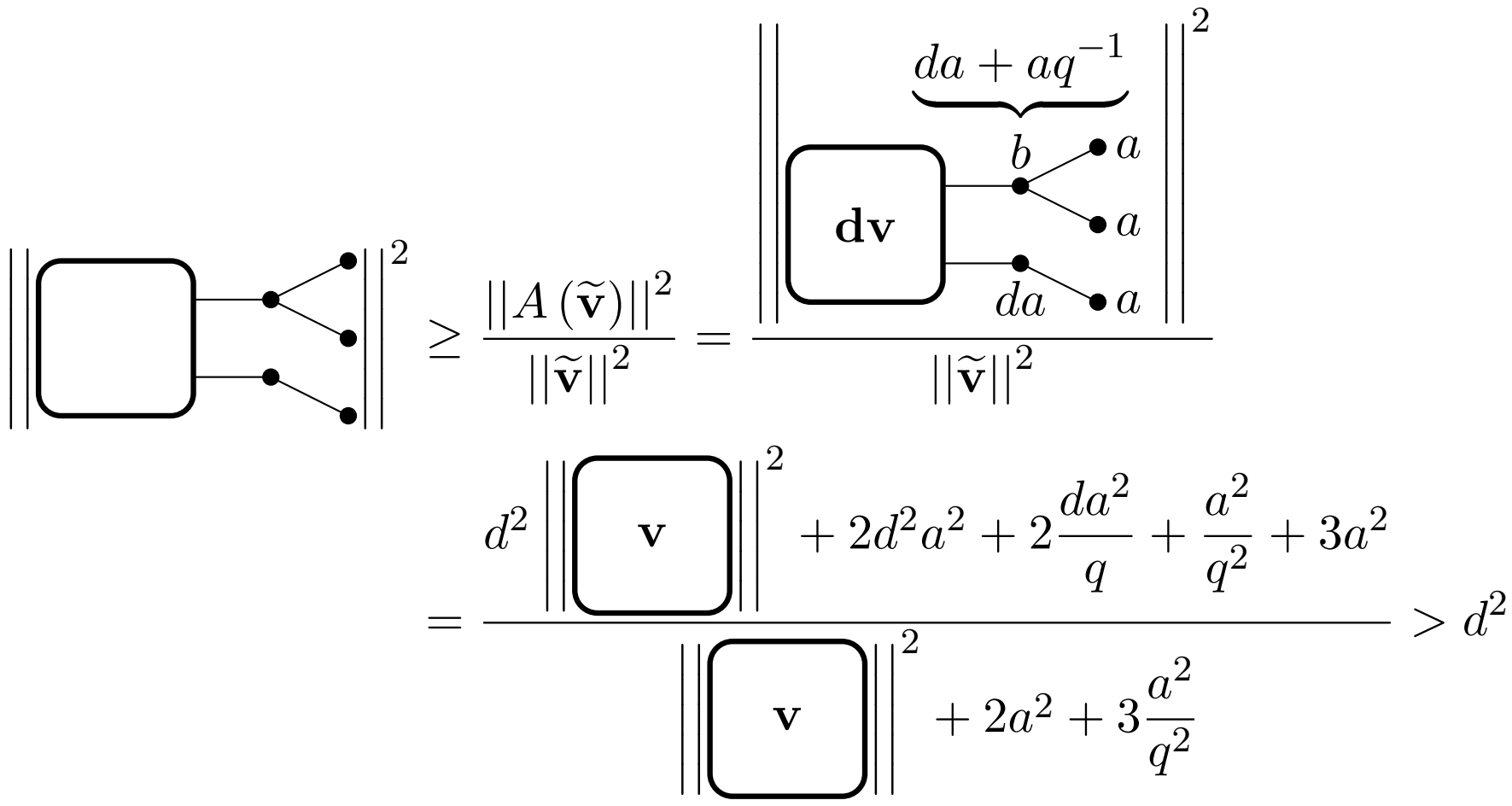

because the inequality

$$d^2a^2 + 2\frac{da^2}{q} + \frac{a^2}{q^2} + 3a^2 > d^2\left(2a^2 + 3\frac{a^2}{q^2}\right)$$

is equivalent to $2q^2 - 3 - 3q^{-2} > 0$. By Theorem A.3, we are finished. □